\pdfoutput=1
\documentclass[11pt]{amsart}
\usepackage[a4paper,margin=29mm]{geometry}
\usepackage[T1]{fontenc}
\usepackage{newpxtext}
\usepackage{mathtools}
\usepackage{newpxmath}
\usepackage{microtype}
\usepackage{xcolor}
\definecolor{InternalRed}{RGB}{0,55,100}
\definecolor{ExternalGreen}{RGB}{0,85,70}
\usepackage[colorlinks=true,linkcolor=InternalRed,citecolor=ExternalGreen,
urlcolor=ExternalGreen]{hyperref}
\usepackage{enumitem}
\usepackage{needspace}
\numberwithin{equation}{section}
\newtheorem{theorem}{Theorem}[section]
\newtheorem{maintheorem}{Theorem}

\newtheorem{proposition}[theorem]{Proposition}
\newtheorem{lemma}[theorem]{Lemma}
\newtheorem{corollary}[theorem]{Corollary}
\theoremstyle{definition}
\newtheorem{definition}[theorem]{Definition}
\newtheorem{conjecture}[theorem]{Conjecture}
\newtheorem{condition}[theorem]{Condition}
\newtheorem{example}[theorem]{Example}
\theoremstyle{remark}

\DeclareMathOperator{\Supp}{Supp}
\DeclareMathOperator{\Spec}{Spec}
\DeclareMathOperator{\Max}{Max}
\DeclareMathOperator{\Min}{Min}
\DeclareMathOperator{\Perf}{Perf}
\DeclareMathOperator{\Loc}{Loc}
\DeclareMathOperator{\Thick}{Thick}
\DeclareMathOperator{\Hom}{Hom}
\DeclareMathOperator{\Ext}{Ext}
\DeclareMathOperator{\Frac}{Frac}
\DeclareMathOperator{\im}{im}
\DeclareMathOperator{\coker}{coker}
\DeclareMathOperator{\cone}{cone}
\DeclareMathOperator{\res}{res}
\newcommand{\RHom}{\mathbf R\!\Hom}
\newcommand{\Lotimes}{\otimes_R^{\mathbf L}}
\newcommand{\D}{\mathbf D}
\newcommand{\fg}{\mathrm{fg}}
\newcommand{\Lex}{\Loc_{\mathrm{ex}}}
\newcommand{\Lten}{\Loc_{\mathrm{ex}}^{\otimes}}
\newcommand{\Lhom}{\Loc_{\mathrm{ex}}^{\Hom}}
\newcommand{\Mod}{\operatorname{Mod}}
\newcommand{\p}{\mathfrak p}
\newcommand{\q}{\mathfrak q}
\newcommand{\m}{\mathfrak m}
\makeatletter
\renewcommand{\@settitle}{%
  \begin{center}
    \normalfont\scshape%
    \fontsize{16}{21}\selectfont
    \@title\par
  \end{center}%
}
\makeatother

\title[A relative telescope-type conjecture]{A Relative Telescope-Type Conjecture:\\[2pt]
Support and Generation in\\[2pt]
Half-Bounded Derived Categories}
\author{Liangyu Zhu}
\address{School of Mathematical Sciences, Anhui University, Hefei 230601, China}
\email{a126112006@stu.ahu.edu.cn}
\date{3 November 2026}
\hypersetup{pdftitle={A Relative Telescope-Type Conjecture: Support and Generation in Half-Bounded Derived Categories},pdfauthor={Liangyu Zhu}}

\begin{document}
\begin{abstract}
Let $R$ be a commutative noetherian ring. We ask whether complexes
with finitely generated cohomology and a fixed support can be
generated by relatively compact objects, using coproducts which
exist in the half-bounded derived category. We prove that
$\D^+_{\fg,W}(R)=\Loc_{\mathrm{ex}}(\D^b_{\fg,W}(R))$ whenever
$W\subseteq\Spec R$ is specialization closed and has finite
dimension. The same equality holds for every specialization-closed
subset of $\Spec B$ when $B$ is a finite-type algebra over a
noetherian domain $A$ of finite character and either $A$ or $B$
is regular. These equalities classify the thick Hom-submodules
closed under existing coproducts. We also prove the corresponding
bounded-above statement from the support results of
Matsui--Takahashi. For an arbitrary noetherian ring, the
bounded-below generation conjecture reduces to torsion complexes
over quotient domains, with their support unchanged. Finally,
over a regular noetherian domain of infinite dimension, we give
a bounded-below object which is generated by bounded objects
with the same support but has no bounded-below representative
by finitely generated modules. The comparison functor from the
bounded-below derived category of finite modules is therefore
neither full nor essentially surjective.
\end{abstract}
\maketitle

\section{Introduction}

Let $R$ be a commutative noetherian ring. Neeman's theorem identifies
the smashing subcategories of $\D(R)$ with the specialization-closed
subsets of $\Spec R$ \cite{Neeman92}. In particular, complexes with
cohomology supported in such a subset $W$ are generated by perfect
complexes with support in $W$, using triangles and arbitrary coproducts
\cite{Neeman92,BIK08}. The perfect complexes are the compact objects
of $\D(R)$; they are also its rigid objects. This is a basic example
of rigid tensor triangular geometry. Balmer--Favi explain how the
telescope conjecture and its associated tensor idempotents fit into
this setting \cite{BF11}.

What happens if we require the cohomology modules to remain finite
and the complexes to remain bounded on one side? Write
$\D^-_{\fg}(R)$ and $\D^+_{\fg}(R)$ for the corresponding full
subcategories of $\D(R)$. The category $\D^-_{\fg}(R)$ is tensor
triangulated, but it contains objects which are not rigid. Its thick
tensor ideals therefore belong to non-rigid tensor triangular
geometry. Matsui--Takahashi studied these ideals in \cite{MT17}.
Even over a discrete valuation ring, ideals with the same ordinary
support can differ considerably; the growth of torsion in cohomology
distinguishes them \cite{SZ26}.

Here we ask a generation question in each of the two half-bounded
categories. We use only coproducts which exist in the category in
question. For a collection $\mathcal E$, let $\Lex(\mathcal E)$
denote its thick closure under these coproducts.

\Needspace{7\baselineskip}
For a
specialization-closed $W\subseteq\Spec R$, does the following
equality hold?
\[
 \D^+_{\fg,W}(R)=\Lex(\D^b_{\fg,W}(R))
\]
Can we build a bounded-below complex from bounded complexes
with the same support, without leaving $\D^+_{\fg}(R)$?
We prove this when $W$ has finite dimension. We also prove it for
certain rings of infinite dimension, including singular rings.

Why not simply take the bounded truncations? Their usual homotopy
colimit recovers $X^{\bullet}$ in $\D(R)$, but the coproduct in the
telescope need not belong to $\D^+_{\fg}(R)$. For example, if
$X^{\bullet}\in\D^{\geq a}(R)$ and $H^i(X^{\bullet})\ne0$, then
\[
 H^i\!\left(\bigoplus_{n\geq a}\tau_{\leq n}X^{\bullet}\right)
 =\bigoplus_{n\geq\max\{a,i\}}H^i(X^{\bullet})
\]
is not finitely generated. The theorem in $\D(R)$ leaves this
particular difficulty unresolved.

One might instead look for a bounded-below representative by finite
modules supported in $W$. Separating its even and odd terms would
then prove the required generation statement
(Lemma~\ref{lem:blocks}). Such representatives exist locally when
the cohomology has finite length. Can we lift these local
representatives and combine them? Can we do so without introducing
infinitely many summands in a fixed degree? These questions lead
to the proof of Theorem~\ref{thm:positive-intro}.

At a minimal prime $\p$ of $\Supp X^{\bullet}$, we lift a local
representative starting in degree $\inf X^{\bullet}_\p$.
The lifting requires two choices of denominators: first for maps
between bounded summands, then for the map from their cone to
$X^{\bullet}$ (Lemma~\ref{lem:scalar-spreading}). The starting
degrees ensure that only finitely many minimal primes contribute
in each degree. We may therefore take their coproduct in
$\D^+_{\fg}(R)$. The cone has smaller-dimensional support, so
induction applies.

There is a second difficulty. A representative by finite modules
need not exist globally, even over a regular noetherian ring.
Theorem~\ref{thm:representability-intro} gives an example for which
the generation statement holds nevertheless. Thus the local
representatives used in the proof cannot always be replaced by
one global representative of the same kind.

The category $\D^+_{\fg}(R)$ is generally not closed under derived
tensor product. We use the Hom action
\begin{equation}\label{eq:action-intro}
 \begin{aligned}
 \odot:\D^-_{\fg}(R)^{\mathrm{op}}\times\D^+_{\fg}(R)
 &\longrightarrow\D^+_{\fg}(R),\\
 A^{\bullet}\odot X^{\bullet}&=\RHom_R(A^{\bullet},X^{\bullet})
 \end{aligned}
\end{equation}
of \cite[Lemma~6.4 and Appendix~B.4]{Zhu26a}. A thick
Hom-submodule is a thick subcategory closed under this action.
Compactness, measured against coproducts existing in the indicated
category, gives
\begin{equation}\label{eq:cores-intro}
 \D^-_{\fg}(R)^c=\Perf(R),\qquad
 \D^+_{\fg}(R)^c=\D^b_{\fg}(R).
\end{equation}
The first equality and this use of compactness are in
\cite[Proposition~2.1]{MT17}; the second is
Lemma~\ref{lem:relative-compact-objects}. In particular, the compact
objects on the bounded-below side can have infinite projective
dimension.

The telescope conjecture suggests asking whether subcategories
closed under the action and under existing coproducts are generated
by these relatively compact objects. Stevenson's relative telescope
conjecture concerns actions on compactly generated categories and
retains the requirement of a smashing Bousfield localization
\cite[Definition~7.1]{Stevenson13}. Our half-bounded categories have
fewer coproducts. Those which exist are also products
(Lemma~\ref{lem:admissible-coproducts}), so right orthogonals are
automatically closed under them. Consequently, relative smashing
in Definition~\ref{def:smashing} amounts to closure under existing
coproducts. The resulting telescope-type question is the generation
question just described; its precise formulation is
Definition~\ref{def:telescope-type}.

On the bounded-above side the answer follows from known support
results and a short argument with complexes. Matsui--Takahashi
classify the thick tensor ideals generated by bounded or perfect
objects \cite[Theorem~2.12 and Corollary~2.14]{MT17}. A supported
representative and its even and odd terms then give
\[
 \D^-_{\fg,W}(R)=\Lex(\D^b_{\fg,W}(R))=\Lten(\Perf_W(R)).
\]
We include the argument in Proposition~\ref{thm:negative}.
The coproduct condition is essential here: this equality says
nothing about the other thick tensor ideals studied in
\cite{MT17,SZ26}.

For the bounded-below category we prove the following.

\begin{maintheorem}
\label{thm:positive-intro}
Let $W\subseteq\Spec R$ be specialization closed. If
$\sup_{\p\in W}\dim(R/\p)<\infty$, then
\[
 \D^+_{\fg,W}(R)=\Lex(\D^b_{\fg,W}(R)).
\]
Let $\mathcal S\subseteq\D^+_{\fg}(R)$ be a thick Hom-submodule closed
under existing coproducts. If
$\dim(R/\p)$ is bounded on $\Supp\mathcal S$, then
\[
 \mathcal S=\D^+_{\fg,\Supp\mathcal S}(R)
 =\Lhom(\mathcal S\cap\D^b_{\fg}(R)).
\]
\end{maintheorem}

The first equality follows by removing minimal support primes as
above. For the second, we prove that
$R/\p\in\mathcal S$ whenever $\p\in\Supp\mathcal S$.
Over the domain $R/\p$, choose free modules which become bases
after passage to its fraction field. Their cone has smaller
support, and the same induction applies
(Theorem~\ref{thm:coheight}).

When $R$ admits a dualizing complex, duality already transfers the
bounded-above generation statement to $\D^+_{\fg}(R)$
\cite[Tag~0A7C]{Stacks}. Theorem~\ref{thm:positive-intro} treats
finite-dimensional support over an arbitrary noetherian ring.
To go beyond finite dimension, we use dualizing complexes after
localization. Their injective dimensions need not have a uniform
bound. The next theorem shows how finite character allows us to
use the finitely many local bounds needed for each bounded torsion
complex.

\begin{maintheorem}\label{thm:infinite-intro}
Let $A$ be a noetherian domain of finite character and let $B$ be a
noetherian $A$-algebra. Suppose that there is
$\omega^{\bullet}\in\D^b_{\fg}(B)$ such that
$\omega^{\bullet}_{\m}$ is dualizing over $B_{\m}$ whenever
$\m\in\Max A$ and $B_{\m}\ne0$. For every specialization-closed
$W\subseteq\Spec B$,
\[
 \D^+_{\fg,W}(B)=\Lex(\D^b_{\fg,W}(B)).
\]
If $W\subseteq\{\q\in\Spec B:\q\cap A\ne(0)\}$, every object of
$\D^+_{\fg,W}(B)$ has a bounded-below representative by finite
$W$-supported modules.
\end{maintheorem}

Finite character means that each nonzero element of $A$ belongs
to only finitely many maximal ideals. If $B$ is a finite-type
$A$-algebra and either $A$ or $B$ is regular, the hypothesis on
$\omega^{\bullet}$ holds (Corollary~\ref{thm:finite-character}).
In these cases, every thick Hom-submodule
$\mathcal S\subseteq\D^+_{\fg}(B)$ closed under existing
coproducts satisfies
\[
 \mathcal S=\D^+_{\fg,\Supp\mathcal S}(B)
 =\Lhom(\mathcal S\cap\D^b_{\fg}(B)).
\]
For Nagata's regular domain $A$ in Example~\ref{ex:nagata}, this
applies, for instance, to the singular ring
\[
 B=A[t,\varepsilon]/(\varepsilon^2),\qquad
 \p=(t,\varepsilon),\qquad \dim(B/\p)=\infty.
\]
This example is not covered by the finite-dimensional induction.

What remains for a general noetherian ring? The same argument
isolates a question about torsion complexes over the domains
$R/\p$.

\begin{maintheorem}
\label{thm:remaining-intro}
Let $\mathcal S\subseteq\D^+_{\fg}(R)$ be a thick Hom-submodule closed
under existing coproducts, and set
$\mathcal U=\Lhom(\mathcal S\cap\D^b_{\fg}(R))$. For every prime $\p$
with $\dim(R/\p)<\infty$,
\[
 \begin{aligned}
 \p\in\Supp\mathcal S&\ \Longrightarrow\ R/\p\in\mathcal S\cap\D^b_{\fg}(R),\\
 \RHom_R(R/\p,X^{\bullet})&\in\mathcal U\qquad(X^{\bullet}\in\mathcal S).
 \end{aligned}
\]
For a fixed $R$, Conjecture~\ref{conj:positive} is equivalent
to the following assertion. For every $\p\in\Spec R$, with $B=R/\p$
and $K=\Frac(B)$,
\[
 \begin{gathered}
 C^{\bullet}\in\D^+_{\fg}(B),\qquad C^{\bullet}\otimes_BK=0\\
 \Longrightarrow\quad
 \res C^{\bullet}\in\Loc_{\mathrm{ex},R}
 \bigl(\D^b_{\fg,\Supp_R(\res C^{\bullet})}(R)\bigr).
 \end{gathered}
\]
\end{maintheorem}

The first two assertions hold for every noetherian ring. The
conjecture remains open in general. In the last assertion it is
important to keep the actual support
$\Supp_R(\res C^{\bullet})$: generation by bounded torsion objects
with larger support would not suffice.

We finish with the promised example concerning representatives
by finite modules.

\begin{maintheorem}
\label{thm:representability-intro}
Let $A$ be the regular noetherian domain in Example~\ref{ex:nagata},
and set $B=A[t]$. There is $C^{\bullet}\in\D^+_{\fg}(B)$ such that
\begin{enumerate}[label=(\roman*),leftmargin=*]
\item $C^{\bullet}$ has a bounded-below representative on which $t$ acts termwise
as zero, and $\Supp_B C^{\bullet}=V(t)$;
\item $C^{\bullet}\in\Lex(\D^b_{\fg,V(t)}(B))$;
\item $C^{\bullet}$ has no bounded-below representative by finite $B$-modules.
\end{enumerate}
For this $B$, the natural functor
\[
 \D^+(B\text{-}\mathrm{mod})\longrightarrow\D^+_{\fg}(B)
\]
is neither full nor essentially surjective.
\end{maintheorem}

Positselski--Schn\"urer study the comparison between derived
categories of finite modules and complexes with finite cohomology
\cite{PS21}. Their quasi-Frobenius counterexample is unbounded
in both directions \cite[Theorem~3.1]{PS21}; for those rings the
bounded-below comparison is an equivalence
\cite[Section~3, (3.2) and remarks~(a)--(b)]{PS21}.
Finite global dimension gives an unbounded equivalence
\cite[Corollaries~4.5 and~5.11]{PS21}.

In Theorem~\ref{thm:representability-intro}, the ring is regular
but has infinite dimension. We choose nonzero extension classes
of increasing degrees at different maximal ideals. A bounded-below
representative by finite modules would give a finite $B$-module
$N$ and a monomorphism $A\to N$ which sends all these extension
classes to zero. By generic freeness this monomorphism splits
as a map of $A$-modules after one localization. Finite character then leaves
a maximal ideal where the splitting and a nonzero extension class
both survive, a contradiction. Theorem~\ref{thm:no-finite-model}
contains the construction. The failure of fullness in
Corollary~\ref{cor:comparison-failure} follows by the usual fiber
argument, as in \cite[proof of Theorem~3.1]{PS21}.

For comparison, Matsui--Takahashi give a classification in
$\D^+(\Mod R)$ under their $n$-uniformity condition, closure under terms of injective
resolutions, and closure under existing sums and products
\cite[Theorem~6.7(1)]{MT22}; their coproduct convention is explicit
in \cite[Convention~2.1]{MT22}. The injective-resolution condition
does not pass to our subcategories of $\D^+_{\fg}(R)$.
Over a discrete valuation ring, the finite-cohomology torsion
subcategory contains the residue field but not its injective hull
(Example~\ref{ex:not-bousfield}). The proofs below instead keep
finite cohomology and the chosen support at every step.

\section{Preliminaries}\label{sec:preliminaries}

Throughout, rings are commutative and noetherian. Subcategories are full
and closed under isomorphisms. Write $\Mod R$ for all $R$-modules and
$R\text{-}\mathrm{mod}$ for the finitely generated ones; for modules,
``finite'' means finitely generated. We denote complexes by
$X^{\bullet},Y^{\bullet},\ldots$, and their degree-$i$ terms by
$X^i,Y^i,\ldots$. We use cohomological indexing:
\[
 H^i(X^{\bullet}[s])=H^{i+s}(X^{\bullet}),\qquad
 \inf X^{\bullet}=\inf\{i:H^i(X^{\bullet})\ne0\},\quad
 \sup X^{\bullet}=\sup\{i:H^i(X^{\bullet})\ne0\},
\]
with $\inf0=+\infty$ and $\sup0=-\infty$.
The subscript $\fg$ requires every cohomology module to be finite;
in particular, $\D_{\fg}(R)$ denotes the full subcategory of the
unbounded derived category with this property. The superscripts $-$, $+$, and $b$ mean bounded above, bounded below,
and bounded, respectively. The functors $\tau_{\leq b}$ and
$\tau_{\geq a}$ are the canonical cohomological truncations.
All morphisms are taken in $\D(R)=\D(\Mod R)$.
We write $\RHom_R$ and $\otimes_R^{\mathbf L}$ for the derived functors.
In arguments involving a quotient $R\to B$, a $B$-complex is also
viewed in $\D(R)$ by restriction of scalars, denoted by $\res$ when
the base ring must be distinguished.

The category $\Perf(R)$ consists of the objects represented by bounded
complexes of finite projective modules. A representative by finite modules may be unbounded; each of its terms is
finitely generated. The natural functor
$\D^-(R\text{-}\mathrm{mod})\to\D^-_{\fg}(R)$ is an equivalence
\cite[Proposition~2.1]{PS21}; its positive analogue need not be
essentially surjective (Theorem~\ref{thm:no-finite-model}).
We set $k(\p)=R_\p/\p R_\p$ and
\[
 \begin{aligned}
 \Supp_RM&=\{\p\in\Spec R:M_\p\ne0\},\\
 \Supp_R X^{\bullet}&=\bigcup_i\Supp_RH^i(X^{\bullet}),\\
 \Supp_R\mathcal E&=\bigcup_{X^{\bullet}\in\mathcal E}\Supp_R X^{\bullet}.
 \end{aligned}
\]
We omit the ring subscript when it is unambiguous. A subset
$W\subseteq\Spec R$ is specialization closed if it is closed under
passage from a prime ideal to a larger prime ideal. For
$\epsilon\in\{-,+,b\}$, set
\[
 \begin{aligned}
 \D^\epsilon_{\fg,W}(R)&=
 \{X^{\bullet}\in\D^\epsilon_{\fg}(R):\Supp_R X^{\bullet}\subseteq W\},\\
 \Perf_W(R)&=\Perf(R)\cap\D^b_{\fg,W}(R).
 \end{aligned}
\]
Such objects are called \emph{$W$-supported}.
The \emph{coheight} of $\p$ is $\dim(R/\p)$, and
\begin{equation}\label{eq:dimension}
 d(W)=\sup_{\p\in W}\dim(R/\p),\qquad d(\varnothing)=-1.
\end{equation}
For a multiplicative subset $\Phi\subseteq R$, write
$R_\Phi=\Phi^{-1}R$, $X^{\bullet}_\Phi=\Phi^{-1}X^{\bullet}$, and
\[
 W_\Phi=\{\p R_\Phi:\p\in W,\ \p\cap\Phi=\varnothing\}.
\]

\begin{definition}\label{def:closures}
Fix $\epsilon\in\{-,+\}$. A family $(X^{\bullet}_\lambda)_{\lambda\in\Lambda}$,
indexed by a set, in
$\D^\epsilon_{\fg}(R)$ is \emph{admissible} if its coproduct exists in
that category. A thick subcategory is \emph{relatively localizing} if it
is closed under coproducts of all admissible families of its objects.
We denote the smallest such subcategory containing $\mathcal E$ by
$\Lex(\mathcal E)$; the ambient category is specified in each statement.

A thick tensor ideal $\mathcal S\subseteq\D^-_{\fg}(R)$ satisfies
$A^{\bullet}\Lotimes X^{\bullet}\in\mathcal S$ for $A^{\bullet}\in\D^-_{\fg}(R)$ and
$X^{\bullet}\in\mathcal S$. A \emph{thick Hom-submodule}
$\mathcal S\subseteq\D^+_{\fg}(R)$ satisfies
\[
 \RHom_R(A^{\bullet},X^{\bullet})\in\mathcal S
 \qquad(A^{\bullet}\in\D^-_{\fg}(R),\ X^{\bullet}\in\mathcal S).
\]
The symbols $\Lten$ and $\Lhom$ denote relatively localizing closure
with the tensor-ideal and Hom-submodule conditions, respectively.
The symbol $\Thick^\otimes$ denotes thick tensor-ideal closure alone.
A ring subscript on $\Loc$ specifies the base ring.
\end{definition}

\begin{definition}\label{def:relative-compact}
An object $P^{\bullet}\in\D^\epsilon_{\fg}(R)$, $\epsilon\in\{-,+\}$, is
\emph{relatively compact} if, for every admissible family, the canonical
map
\[
 \bigoplus_\lambda\Hom_{\D(R)}(P^{\bullet},X^{\bullet}_\lambda)
 \longrightarrow
 \Hom_{\D(R)}\!\left(P^{\bullet},\bigoplus_\lambda X^{\bullet}_\lambda\right)
\]
is an isomorphism. Write $\D^\epsilon_{\fg}(R)^c$ for these objects.
A relatively localizing thick tensor ideal
$\mathcal S\subseteq\D^-_{\fg}(R)$ is \emph{relatively compactly
generated} if
$\mathcal S=\Lten(\mathcal S\cap\D^-_{\fg}(R)^c)$.
For a relatively localizing thick Hom-submodule of $\D^+_{\fg}(R)$,
the corresponding requirement is
$\mathcal S=\Lhom(\mathcal S\cap\D^+_{\fg}(R)^c)$.
\end{definition}

\Needspace{10\baselineskip}
For coproducts, we follow the convention of \cite[Section~2]{Zhu26a}.

\begin{lemma}\label{lem:admissible-coproducts}
For a family $(X^{\bullet}_\lambda)_{\lambda\in\Lambda}$ in
$\D^\epsilon_{\fg}(R)$, $\epsilon\in\{-,+\}$, the following are equivalent:
\begin{enumerate}[label=(\roman*),leftmargin=*]
\item its coproduct exists in $\D^\epsilon_{\fg}(R)$;
\item its product exists in $\D^\epsilon_{\fg}(R)$;
\item the family has a common cohomological bound of the indicated sign, and
\begin{equation}\label{eq:admissible}
 \#\{\lambda:H^i(X^{\bullet}_\lambda)\ne0\}<\infty
 \qquad(i\in\mathbb Z).
\end{equation}
\end{enumerate}
When these conditions hold, the internal coproduct and product agree
with those in $\D(R)$, and the canonical map from the coproduct to the
product is an isomorphism.
\end{lemma}

\begin{lemma}\label{lem:relative-compact-objects}
The relatively compact objects are
\[
 \D^-_{\fg}(R)^c=\Perf(R),\qquad
 \D^+_{\fg}(R)^c=\D^b_{\fg}(R).
\]
\end{lemma}
\begin{proof}
The equality for $\D^-_{\fg}(R)$ is \cite[Proposition~2.1]{MT17}.
The zero object is compact. Let $0\ne C^{\bullet}\in\D^b_{\fg}(R)$
and $b=\sup C^{\bullet}$.
For an admissible family in $\D^+_{\fg}(R)$, Lemma~\ref{lem:admissible-coproducts}
implies that only finitely many indices satisfy
$\tau_{\leq b}X^{\bullet}_\lambda\ne0$.
Since $\Hom(C^{\bullet},Y^{\bullet})=
\Hom(C^{\bullet},\tau_{\leq b}Y^{\bullet})$, the compactness map reduces
to a finite direct sum and is an isomorphism.

If a relatively compact $C^{\bullet}\in\D^+_{\fg}(R)$ were unbounded
above, the family $X^{\bullet}_n=\tau_{\geq n}C^{\bullet}$, $n\geq0$,
would be admissible, and every canonical map
$q_n:C^{\bullet}\to X^{\bullet}_n$ would be nonzero.
The product property gives a map with components $(q_n)_n$, whereas
compactness and Lemma~\ref{lem:admissible-coproducts} give
\[
 \Hom\!\left(C^{\bullet},\prod_{n\geq0}X^{\bullet}_n\right)
 =\Hom\!\left(C^{\bullet},\bigoplus_{n\geq0}X^{\bullet}_n\right)
 \cong\bigoplus_{n\geq0}\Hom(C^{\bullet},X^{\bullet}_n).
\]
Every map on the right has only finitely many nonzero components,
a contradiction.
\end{proof}

We recall the following properties of derived Hom from
\cite[Lemma~6.4 and Appendix~B.4]{Zhu26a}.

\begin{lemma}\label{lem:hom-localization}
For $A^{\bullet}\in\D^-_{\fg}(R)$,
$X^{\bullet}\in\D^+_{\fg}(R)$ and a multiplicative subset $\Phi$,
\begin{equation}\label{eq:localization}
 \Phi^{-1}\RHom_R(A^{\bullet},X^{\bullet})\simeq
 \RHom_{R_\Phi}(A^{\bullet}_\Phi,X^{\bullet}_\Phi).
\end{equation}
In particular, for $E^{\bullet}\in\D^b_{\fg}(R)$,
\begin{equation}\label{eq:bounded-hom-localization}
 \Hom_{\D(R)}(E^{\bullet},X^{\bullet})_\Phi\cong
 \Hom_{\D(R_\Phi)}(E^{\bullet}_\Phi,X^{\bullet}_\Phi).
\end{equation}
Tensor product on $\D^-_{\fg}(R)$ and the action
$A^{\bullet}\odot X^{\bullet}=\RHom_R(A^{\bullet},X^{\bullet})$
preserve admissible coproducts in each variable, with a coproduct in
the contravariant variable taken to the corresponding product.
\end{lemma}

The next equivalence follows from
\cite[Fact~3.5 and Proposition~3.6]{SW17} and Nakayama's lemma.

\begin{lemma}\label{lem:residue-detection}
For $X^{\bullet}\in\D_{\fg}(R)$ and $\p\in\Spec R$,
\begin{equation}\label{eq:detection}
 \RHom_{R_\p}(k(\p),X^{\bullet}_\p)=0
 \quad\Longleftrightarrow\quad X^{\bullet}_\p=0.
\end{equation}
\end{lemma}

\begin{definition}\label{def:smashing}
Fix $\epsilon\in\{-,+\}$ and a thick subcategory
$\mathcal S\subseteq\D^\epsilon_{\fg}(R)$. Its \emph{ordinary right
orthogonal} is
\[
 \mathcal S^\perp=
 \{Y^{\bullet}\in\D^\epsilon_{\fg}(R):
 \Hom_{\D(R)}(X^{\bullet},Y^{\bullet}[n])=0\text{ for all }X^{\bullet}\in\mathcal S,\ n\in\mathbb Z\}.
\]
A thick tensor ideal, or thick Hom-submodule, is \emph{relative smashing}
if both $\mathcal S$ and $\mathcal S^\perp$ are closed under existing
coproducts in the specified ambient category.
\end{definition}

\begin{proposition}\label{prop:orthogonal}
For either sign $\epsilon$, the ordinary right orthogonal of any collection in
$\D^\epsilon_{\fg}(R)$ is closed under admissible coproducts.
A thick tensor ideal or thick Hom-submodule is therefore relative
smashing if and only if it is relatively localizing.
\end{proposition}
\begin{proof}
For an admissible family $(Y^{\bullet}_\lambda)$,
\[
 \begin{aligned}
 \Hom_{\D(R)}\!\left(X^{\bullet},\bigoplus_\lambda Y^{\bullet}_\lambda[n]\right)
 &\cong\Hom_{\D(R)}\!\left(X^{\bullet},\prod_\lambda Y^{\bullet}_\lambda[n]\right)\\
 &\cong\prod_\lambda\Hom_{\D(R)}(X^{\bullet},Y^{\bullet}_\lambda[n]).
 \end{aligned}
\]
\end{proof}

\begin{definition}[Relative telescope-type property]\label{def:telescope-type}
The category $\D^-_{\fg}(R)$ satisfies the \emph{relative telescope-type
property} if every relative smashing thick tensor ideal is relatively
compactly generated. The category $\D^+_{\fg}(R)$ satisfies this property
if every relative smashing thick Hom-submodule is relatively compactly
generated.
\end{definition}

\begin{example}\label{ex:not-bousfield}
Let $(R,\m)$ be a discrete valuation ring and
$\mathcal S=\D^\epsilon_{\fg,\{\m\}}(R)$, for either sign.
Then $\mathcal S$ is a proper relative smashing subcategory, but
$\mathcal S^\perp=0$ by Lemma~\ref{lem:residue-detection}:
every nonzero object of $\D^\epsilon_{\fg}(R)$ has $\m$ in its support.
A Bousfield localization with these local objects would be zero and
would have kernel the whole category.

For $\epsilon=+$, this example also distinguishes the present closure
conditions from those in \cite[Theorem~6.7]{MT22}. The subcategory
contains $R/\m$, but its injective hull
\[
 E_R(R/\m)\cong\Frac(R)/R
\]
is not finitely generated. Viewed inside $\D^+(\Mod R)$, the
subcategory therefore fails the injective-resolution closure assumed
there.
\end{example}

\begin{conjecture}[Relative telescope-type conjecture in the bounded-below category]\label{conj:positive}
Every relatively localizing thick Hom-submodule
$\mathcal S\subseteq\D^+_{\fg}(R)$ satisfies
\[
 \mathcal S=\Lhom(\mathcal S\cap\D^b_{\fg}(R)).
\]
\end{conjecture}

\section{Generation in the bounded-above category}

For a specialization-closed $W$, set
\[
 \Mod_W R=\{M:\Supp_R M\subseteq W\},\qquad
 \D_W(\Mod R)=\{X^{\bullet}\in\D(R):\Supp_R X^{\bullet}\subseteq W\}.
\]
The torsion-derived comparison is
\begin{equation}\label{eq:torsion-comparison}
 \D(\Mod_W R)\xrightarrow{\sim}\D_W(\Mod R),
\end{equation}
as in the proof of \cite[Theorem~9.1]{BIK08}. It applies to arbitrary
commutative noetherian rings. For closed support, see also
\cite[Lemma~47.10.1, Tag~0955]{Stacks}.

\begin{lemma}\label{lem:negative-model}
Let $W\subseteq\Spec R$ be specialization closed.
If $X^{\bullet}\in\D^-_{\fg,W}(R)$ and $H^i(X^{\bullet})=0$ for $i>b$, then $X^{\bullet}$ has a representative
$Q^{\bullet}$ with $Q^i=0$ for $i>b$ and with all $Q^i$ finitely generated and
$W$-supported. If $X^{\bullet}$ is bounded, such a representative can be chosen
bounded.
\end{lemma}
\begin{proof}
By \eqref{eq:torsion-comparison}, choose a complex $C^{\bullet}$ of $W$-supported
modules representing $X^{\bullet}$. Replace it by $\tau_{\leq b}C^{\bullet}$ in
$\Mod_W R$. We construct a subcomplex $Q^{\bullet}\subseteq C^{\bullet}$, with finite terms,
starting in degree $b$ and proceeding downwards. At each step we must
represent the cohomology in the new degree and remove the kernel left
in the preceding degree. Set $Q^{b+1}=0$ and suppose that $Q^{n+1}$
has been chosen.
Let $Q^n$ be generated by cycle representatives of a finite
generating set of $H^n(C^{\bullet})$ and preimages of a finite generating set of
\[
 Q^{n+1}\cap B^{n+1}(C^{\bullet}).
\]
The latter module is finite by noetherianity. Then $d(Q^n)\subseteq Q^{n+1}$.
The chosen cycles make $H^n(Q^{\bullet})\to H^n(C^{\bullet})$
surjective. At the next step we arrange
\[
 Q^n\cap B^n(C^{\bullet})\subseteq d(Q^{n-1}),
\]
so the same map is injective. Hence $Q^{\bullet}\to C^{\bullet}$
is a quasi-isomorphism. Each $Q^n$ is a finite submodule of the
$W$-supported module $C^n$.
For bounded $X^{\bullet}$, a good lower truncation of $Q^{\bullet}$ is bounded and retains
these properties.
\end{proof}

\begin{lemma}\label{lem:blocks}
Let $W\subseteq\Spec R$ be specialization closed.
Fix $\epsilon\in\{-,+\}$, and let $Q^{\bullet}$ be a complex of finite
$W$-supported modules, bounded above for $\epsilon=-$ and bounded below
for $\epsilon=+$.
Then its image in $\D^\epsilon_{\fg}(R)$ belongs to $\Lex(\D^b_{\fg,W}(R))$.
\end{lemma}
\begin{proof}
We separate the even and odd degrees. Fix an integer $r$ and write
$I=r+2\mathbb Z$. For $i\in I$, define
\[
 G^i=Q^i,\qquad G^{i+1}=\operatorname{im}d_Q^i,
 \qquad d_G^i=d_Q^i,\quad d_G^{i+1}=0.
\]
The inclusions $G^j\subseteq Q^j$ form a subcomplex, since
$d_Q^{i+1}(\operatorname{im}d_Q^i)=0$. Thus
\[
 0\longrightarrow G^\bullet\longrightarrow Q^\bullet
   \longrightarrow B^\bullet\longrightarrow0
\]
is degreewise exact, with
\[
 B^i=0,\qquad B^{i+1}=Q^{i+1}/\operatorname{im}d_Q^i
 \quad(i\in I),\qquad d_B=0.
\]
Each two-term block $Q^i\to\operatorname{im}d_Q^i$ has only
$\ker d_Q^i$ as cohomology, in degree $i$. Consequently
\begin{equation}\label{eq:blocks}
\begin{aligned}
 G^\bullet&\simeq\bigoplus_{i\in I}(\ker d_Q^i)[-i],\\
 B^\bullet&=\bigoplus_{i\in I}
       (Q^{i+1}/\operatorname{im}d_Q^i)[-i-1].
\end{aligned}
\end{equation}
The stalks in \eqref{eq:blocks} have the same one-sided bound as
$Q^\bullet$, and in each degree there is at most one summand. Both
coproducts are therefore admissible. Their summands lie in
$\D^b_{\fg,W}(R)$, so the triangle associated to the short exact
sequence gives the assertion.
\end{proof}

\begin{proposition}\label{thm:negative}
For every specialization-closed $W\subseteq\Spec R$,
\begin{equation}\label{eq:negative-generation}
 \D^-_{\fg,W}(R)=\Lex(\D^b_{\fg,W}(R))=\Lten(\Perf_W(R)).
\end{equation}
For every relatively localizing thick tensor ideal $\mathcal S\subseteq\D^-_{\fg}(R)$,
\begin{equation}\label{eq:negative-classification}
 \mathcal S=\D^-_{\fg,\Supp\mathcal S}(R)
 =\Lten(\mathcal S\cap\Perf(R)).
\end{equation}
In particular, for any collection $\mathcal E\subseteq\D^-_{\fg}(R)$,
\begin{equation}\label{eq:negative-arbitrary}
 \Lten(\mathcal E)=\D^-_{\fg,\Supp\mathcal E}(R).
\end{equation}
\end{proposition}
\begin{proof}
The first equality of \eqref{eq:negative-generation} follows from
Lemmas~\ref{lem:negative-model} and~\ref{lem:blocks}. By
\cite[Corollary~2.14]{MT17},
\[
 \D^b_{\fg,W}(R)\subseteq\Thick^{\otimes}_{\D^-_{\fg}(R)}(\Perf_W(R)).
\]
This gives the second equality; the reverse support containment is preserved
by shifts, cones, summands, tensor products, and admissible coproducts.

Set $W=\Supp\mathcal S$. For every $\p\in W$, choose
$X^{\bullet}\in\mathcal S$ with $\p\in\Supp X^{\bullet}$. Then
$R/\p\in\Thick^\otimes_{\D^-_{\fg}(R)}(X^{\bullet})\subseteq\mathcal S$
by \cite[Corollary~2.11]{MT17}.
Prime filtrations and successive cohomological truncations give
$\D^b_{\fg,W}(R)\subseteq\mathcal S$.
Equation~\eqref{eq:negative-generation} yields
$\D^-_{\fg,W}(R)\subseteq\mathcal S$; the opposite inclusion follows from
the definition of $W$.
Moreover $\Perf_W(R)=\mathcal S\cap\Perf(R)$, proving
\eqref{eq:negative-classification}. Finally the support of
$\Lten(\mathcal E)$ equals $\Supp\mathcal E$, which proves
\eqref{eq:negative-arbitrary}.
\end{proof}

\section{Lifting complexes after localization}

\begin{lemma}\label{lem:local-length-model}
Let $(A,\m)$ be noetherian local. Suppose $Z^{\bullet}\in\D^+_{\fg}(A)$,
$H^i(Z^{\bullet})=0$ for $i<a$, and all $H^i(Z^{\bullet})$ have finite length. There is a
representative $Q^{\bullet}$ of $Z^{\bullet}$ with finite-length terms and $Q^i=0$ for $i<a$.
There is a triangle
\begin{equation}\label{eq:length-triangle}
 F_A^{\bullet}\longrightarrow G_A^{\bullet}\longrightarrow Z^{\bullet}\longrightarrow F_A^{\bullet}[1]
\end{equation}
in which $F_A^{\bullet}$ and $G_A^{\bullet}$ are admissible sums of shifted finite-length
modules, with
\begin{equation}\label{eq:length-bounds}
 \inf G_A^{\bullet}\geq a,\qquad \inf F_A^{\bullet}\geq a+2.
\end{equation}
\end{lemma}
\begin{proof}
Let $E=E_A(A/\m)$ and $D_E=\Hom_A(-,E)$, with the usual reversal of complex
degrees. This is exact. Its restriction to finite-length modules is a
duality, and $M\to D_E^2M$ is an isomorphism for every finite-length $M$.
Indeed $M$ is annihilated by a power of $\m$, and finite-length duality over
the resulting artinian quotient applies
\cite[Section~47.6 and Lemma~47.7.1, Tags~08YW and~08Z2]{Stacks}.

Thus $D_E Z^{\bullet}\in\D^-_{\fg}(A)$ has finite-length cohomology and supremum at
most $-a$. Lemma~\ref{lem:negative-model}, with $W=\{\m\}$, gives a
finite-length term representative $P^{\bullet}$, zero above $-a$. Then
\[
 Q^{\bullet}=D_E P^{\bullet}\simeq D_E^2Z^{\bullet}\simeq Z^{\bullet}.
\]
In Lemma~\ref{lem:blocks}, take $I=a+2\mathbb Z$. Its subcomplex
$G_A^{\bullet}$ and quotient $B_A^{\bullet}$ satisfy
\[
 \begin{aligned}
 G_A^{\bullet}&\simeq\bigoplus_{j\geq0}(\ker d_{Q^{\bullet}}^{a+2j})[-a-2j],\\
 B_A^{\bullet}&=\bigoplus_{j\geq0}
 (Q^{a+2j+1}/\im d_{Q^{\bullet}}^{a+2j})[-a-2j-1].
 \end{aligned}
\]
Rotate the triangle and set $F_A^{\bullet}=B_A^{\bullet}[-1]$. This gives
\eqref{eq:length-triangle} and~\eqref{eq:length-bounds}.
\end{proof}

\begin{lemma}\label{lem:supported-descent}
Let $\Phi\subseteq R$ be multiplicative and let $W\subseteq\Spec R$
be specialization closed. Every finite $R_\Phi$-module $L$ supported
in $W_\Phi$ is the localization of a finite $W$-supported $R$-module.
\end{lemma}
\begin{proof}
Assume $L\ne0$. Contract its minimal support primes to
$\p_1,\ldots,\p_s\in W$, and set $J=\bigcap_i\p_i$.
For some $e>0$, $J_\Phi^eL=0$. Choose a finite $R$-module $L_0$
with $(L_0)_\Phi\cong L$, by descending a finite presentation.
Then
\[
 N=L_0/J^eL_0,\qquad N_\Phi\cong L,\qquad
 \Supp_RN\subseteq V(J)=\bigcup_iV(\p_i)\subseteq W.\qedhere
\]
\end{proof}

\begin{lemma}\label{lem:scalar-spreading}
Let $\Phi\subseteq R$ be multiplicative, let $W\subseteq\Spec R$ be
specialization closed, and let $X^{\bullet}\in\D^+_{\fg}(R)$.
Suppose that $X_\Phi^{\bullet}$ has a representative $Q^{\bullet}$, zero below an integer
$a$, whose terms are finite and $W_\Phi$-supported.
There are $Y^{\bullet}\in\Lex(\D^b_{\fg,W}(R))$ and $h:Y^{\bullet}\to X^{\bullet}$ such that
\begin{equation}\label{eq:scalar-spreading}
 \inf Y^{\bullet}\geq a,\qquad h_\Phi\text{ is an isomorphism}.
\end{equation}
\end{lemma}
\begin{proof}
We will lift the decomposition into even and odd degrees. We must lift it with the prescribed support and lower
degree bound. Since there are infinitely many summands, we cannot
expect a common denominator. We instead clear denominators separately
on each bounded summand.
Apply Lemma~\ref{lem:blocks} to $Q^{\bullet}$, taking $I=a+2\mathbb Z$,
and rotate the resulting triangle. We obtain
\[
 F_\Phi^{0,\bullet}\xrightarrow{f^0}G_\Phi^{0,\bullet}
 \xrightarrow{g^0}X_\Phi^{\bullet}\longrightarrow F_\Phi^{0,\bullet}[1],
 \qquad \inf F_\Phi^{0,\bullet}\geq a+2,\quad \inf G_\Phi^{0,\bullet}\geq a,
\]
where $F_\Phi^{0,\bullet}$ and $G_\Phi^{0,\bullet}$ are admissible sums of shifted finite
modules. Lemma~\ref{lem:supported-descent} lifts their summands to
$W$-supported shifted modules $F_j^{\bullet}$ and $G_i^{\bullet}$. Retaining the shifts, form
\[
 F^{\bullet}=\bigoplus_j F_j^{\bullet},\qquad G^{\bullet}=\bigoplus_i G_i^{\bullet}.
\]
These coproducts are admissible, with the same lower bounds as before. Fix identifications
$F_\Phi^{\bullet}\cong F_\Phi^{0,\bullet}$ and $G_\Phi^{\bullet}\cong G_\Phi^{0,\bullet}$.

For each $i$, use \eqref{eq:bounded-hom-localization} to lift a
scalar multiple of $g^0|_{(G_i^{\bullet})_\Phi}$. Together these maps define
$g:G^{\bullet}\to X^{\bullet}$ with $g_\Phi=g^0\beta$, where $\beta$ is a
diagonal automorphism of $G_\Phi^{\bullet}$ with entries in $\Phi$.
Next, for each $j$, lift a scalar multiple of the component of
$\beta^{-1}f^0$ to a map $f_j:F_j^{\bullet}\to G^{\bullet}$.
Since $(gf_j)_\Phi=0$, the same localization formula supplies
$w_j\in\Phi$ with $w_jgf_j=0$. Replacing $f_j$ by $w_jf_j$ gives
\begin{equation}\label{eq:gf}
 f:F^{\bullet}\longrightarrow G^{\bullet},\qquad gf=0,
 \qquad f_\Phi=\beta^{-1}f^0\gamma,
\end{equation}
where $\gamma$ is a diagonal automorphism of $F_\Phi^{\bullet}$.

Complete $f$ to a triangle
\[
 F^{\bullet}\xrightarrow{f}G^{\bullet}\longrightarrow Y^{\bullet}\xrightarrow{\delta}F^{\bullet}[1].
\]
Since $gf=0$, there is a map $h:Y^{\bullet}\to X^{\bullet}$
extending $g$. This alone does not make $h_\Phi$ an isomorphism.
The isomorphic triangles determined by $f^0$ and $f_\Phi$ give an
isomorphism $t:Y_\Phi^{\bullet}\to X_\Phi^{\bullet}$ extending $g_\Phi$.
We compare $h_\Phi$ with this map. Exactness of Hom gives
\begin{equation}\label{eq:local-correction}
 t-h_\Phi=d^0\delta_\Phi
 \qquad\text{for some }d^0:F_\Phi^{\bullet}[1]\to X_\Phi^{\bullet}.
\end{equation}
Again, clear the denominator separately on each bounded summand
$F_j^{\bullet}[1]$ of the source of $d^0$. There are $s_j\in\Phi$ and $d:F^{\bullet}[1]\to X^{\bullet}$ such that, for
$s=\bigoplus_js_j:F^{\bullet}\to F^{\bullet}$,
\[
 d_\Phi=d^0s_\Phi[1].
\]
Let $Y^{\prime\bullet}=\cone(fs)$ with connecting map $\delta':Y^{\prime\bullet}\to F^{\bullet}[1]$.
The maps $(s,1_{G^{\bullet}})$ extend to a morphism of triangles with third map
$k:Y^{\prime\bullet}\to Y^{\bullet}$, satisfying
\[
 \delta k=s[1]\delta',\qquad k_\Phi\text{ is an isomorphism}.
\]
For $h'=hk+d\delta'$, equation~\eqref{eq:local-correction} yields
\[
 h'_\Phi=h_\Phi k_\Phi+d^0s_\Phi[1]\delta'_\Phi
 =(h_\Phi+d^0\delta_\Phi)k_\Phi=tk_\Phi.
\]
Thus $h'_\Phi$ is an isomorphism. The triangle
$F^{\bullet}\xrightarrow{fs}G^{\bullet}\to Y^{\prime\bullet}\to F^{\bullet}[1]$ gives
\[
 Y^{\prime\bullet}\in\Lex(\D^b_{\fg,W}(R)),\qquad
 \inf Y^{\prime\bullet}\geq\min\{\inf G^{\bullet},\inf F^{\bullet}-1\}\geq a.
\]
The pair $(Y^{\prime\bullet},h')$ therefore has the required properties.
\end{proof}

\begin{lemma}\label{lem:spreading}
Let $X^{\bullet}\in\D^+_{\fg}(R)$ and $\p\in\Min\Supp X^{\bullet}$. Set
$a_\p=\inf X_\p^{\bullet}$. There are
$Y^{\bullet}\in\Lex(\D^b_{\fg,V(\p)}(R))$ and $h:Y^{\bullet}\to X^{\bullet}$ with
\begin{equation}\label{eq:spreading}
 \Supp Y^{\bullet}\subseteq V(\p),\qquad \inf Y^{\bullet}\geq a_\p,
 \qquad h_\p\text{ is an isomorphism}.
\end{equation}
\end{lemma}
\begin{proof}
Every $H^i(X_\p^{\bullet})$ has finite length over $R_\p$.
Lemma~\ref{lem:local-length-model} gives a finite-length representative
of $X_\p^{\bullet}$, zero below $a_\p$. Apply
Lemma~\ref{lem:scalar-spreading} with $\Phi=R\setminus\p$ and
$W=V(\p)$.
\end{proof}

\section{Support of finite dimension}

\begin{theorem}\label{thm:geometric}
If $W$ is specialization closed and $d(W)<\infty$, then
\begin{equation}\label{eq:positive-geometric}
 \D^+_{\fg,W}(R)=\Lex(\D^b_{\fg,W}(R)).
\end{equation}
\end{theorem}
\begin{proof}
Set $\mathcal L=\Lex(\D^b_{\fg,W}(R))$. The inclusion
$\mathcal L\subseteq\D^+_{\fg,W}(R)$ follows from preservation of support.
For the reverse inclusion, we induct on $d(W)$, starting with
$W=\varnothing$.

Let $0\ne X^{\bullet}\in\D^+_{\fg,W}(R)$ and set
$W_{X^{\bullet}}=\Supp X^{\bullet}$. We will remove all minimal primes of
$W_{X^{\bullet}}$ at once. There may be infinitely many, so we first check
that the required coproduct exists in $\D^+_{\fg}(R)$.
For each $\p\in\Min W_{X^{\bullet}}$, Lemma~\ref{lem:spreading} gives a map
$Y^{\bullet}(\p)\to X^{\bullet}$ with lower bound
$a_\p=\inf X_\p^{\bullet}$ and support in $V(\p)\subseteq W$.

Set $a=\inf X^{\bullet}$. For every $n\geq a$,
\begin{equation}\label{eq:min-finite}
 \{\p\in\Min W_{X^{\bullet}}:a_\p\leq n\}
 \subseteq\bigcup_{i=a}^n\Min\Supp_R H^i(X^{\bullet}).
\end{equation}
Indeed, $a_\p\le n$ places $\p$ in the support of some
$H^i(X^{\bullet})$, with $a\le i\le n$. Its minimality in
$W_{X^{\bullet}}$ makes it minimal in that support. The right side of
\eqref{eq:min-finite} is finite. Since
$\inf Y^\bullet(\p)\ge a_\p\ge a$, for each fixed $n$ the ambient
coproduct satisfies
\[
 H^n\!\left(\bigoplus_{\p\in\Min W_{X^\bullet}}
                        Y^\bullet(\p)\right)
   =\bigoplus_{\substack{\p\in\Min W_{X^\bullet}\\a_\p\le n}}
                        H^n(Y^\bullet(\p)),
\]
with only finitely many possibly nonzero terms. Hence the coproduct
is admissible and belongs to $\mathcal L$. Form the triangle
\begin{equation}\label{eq:remove-min}
 \bigoplus_{\p\in\Min W_{X^{\bullet}}}Y^{\bullet}(\p)\longrightarrow X^{\bullet}\longrightarrow C^{\bullet}\longrightarrow.
\end{equation}
Distinct minimal primes of $W_{X^{\bullet}}$ are incomparable. Localizing at a minimal
$\p$ kills all summands except $Y^{\bullet}(\p)$, whose map is an isomorphism. Thus
\begin{equation}\label{eq:cone-support}
 \Supp C^{\bullet}\subseteq W_{X^{\bullet}}\setminus\Min W_{X^{\bullet}}.
\end{equation}

For $\q\in W_{X^\bullet}$, the nonempty set
\[
 \{\p\in W_{X^\bullet}:\p\subseteq\q\}
\]
has a minimal member: strictly descending chains in this set have
length at most $\dim R_\q<\infty$. Such a member is also minimal in
$W_{X^\bullet}$. If $\q\in\Supp C^\bullet$, equation
\eqref{eq:cone-support} makes this inclusion strict. Concatenating
$\p\subsetneq\q$ with a chain above $\q$ gives
\[
 \dim(R/\q)+1\le\dim(R/\p)\le d(W).
\]
Hence
\[
 d(\Supp C^\bullet)\le d(W)-1,\qquad
 \Lex(\D^b_{\fg,\Supp C^\bullet}(R))\subseteq\mathcal L.
\]
If $d(W)>0$, induction gives $C^\bullet\in\mathcal L$.
If $d(W)=0$, its support is empty, so $C^\bullet=0$.
Both the first and third terms of \eqref{eq:remove-min} lie in
$\mathcal L$. Thus $X^{\bullet}\in\mathcal L$.
\end{proof}

\begin{lemma}\label{lem:prime-filtration}
Let $\mathcal S\subseteq\D(R)$ be thick and let $Z\subseteq\Spec R$
be specialization closed.
If $R/\q\in\mathcal S$ for every $\q\in Z$, then
$\D^b_{\fg,Z}(R)\subseteq\mathcal S$.
\end{lemma}
\begin{proof}
A finite module supported in $Z$ has a finite prime filtration with factors
$R/\q$, $\q\in Z$. Successive cohomological truncations give the assertion for bounded objects.
\end{proof}

\begin{theorem}\label{thm:coheight}
Let $\mathcal S\subseteq\D^+_{\fg}(R)$ be a relatively localizing thick
Hom-submodule and set
$\mathcal U=\Lhom(\mathcal S\cap\D^b_{\fg}(R))$. If $\dim(R/\p)<\infty$, then
\begin{align}
 \p\in\Supp\mathcal S&\ \Longrightarrow\ R/\p\in\mathcal S\cap\D^b_{\fg}(R),
 \label{eq:capture}\\
 \RHom_R(R/\p,X^{\bullet})&\in\mathcal U\qquad(X^{\bullet}\in\mathcal S).
 \label{eq:coinduce}
\end{align}
\end{theorem}
\begin{proof}
We prove the two assertions together by induction on
$d=\dim(R/\p)$. The reason for keeping them together is that, after
passing to $B=R/\p$, we will need both $B\in\mathcal S$ and generation
of $\RHom_R(B,X^{\bullet})$ by bounded objects of $\mathcal S$.
Write $Y^{\bullet}=\RHom_R(B,X^{\bullet})$. A bounded-below injective
representative for $X^{\bullet}$ makes $Y^{\bullet}$ a complex of $B$-modules. Its
cohomology is finite over $B$, so $Y^{\bullet}\in\D^+_{\fg}(B)$. Equations
\eqref{eq:localization} and~\eqref{eq:detection} give
\begin{equation}\label{eq:Ysupport}
 \Supp_R Y^{\bullet}\subseteq V(\p)\cap\Supp X^{\bullet},
 \qquad Y_\p^{\bullet}\ne0\ \Longleftrightarrow\ X_\p^{\bullet}\ne0.
\end{equation}

For $d=0$, $B$ is a field. In $\D(B)$,
\[
 Y^{\bullet}\simeq\bigoplus_nH^n(Y^{\bullet})[-n].
\]
This is an admissible sum of finite-dimensional $B$-vector spaces. If
$X_\p^{\bullet}\ne0$, a nonzero summand makes a shift of $B$ a retract of $Y^{\bullet}$.
Thus $B\in\mathcal S\cap\D^b_{\fg}(R)$, and $Y^{\bullet}\in\mathcal U$.
If $X_\p^{\bullet}=0$, then $Y^{\bullet}=0$. Whenever
$\p\in\Supp\mathcal S$, choose $X^{\bullet}\in\mathcal S$ with
$X_\p^{\bullet}\ne0$ to obtain \eqref{eq:capture}.

Suppose $d>0$. Each $\q\supsetneq\p$ satisfies
$\dim(R/\q)\leq d-1$. If $X_\p^{\bullet}=0$, then
$Z=\Supp_R Y^{\bullet}\subseteq V(\p)\setminus\{\p\}$ and $d(Z)\leq d-1$.
Every $\q\in Z$ lies in $\Supp\mathcal S$, so induction and
Lemma~\ref{lem:prime-filtration} give
\[
 \D^b_{\fg,Z}(R)\subseteq\mathcal S\cap\D^b_{\fg}(R).
\]
Theorem~\ref{thm:geometric} gives $Y^{\bullet}\in\mathcal U$.

It remains to consider $X_\p^{\bullet}\ne0$. In this case
$V(\p)\subseteq\Supp X^{\bullet}\subseteq\Supp\mathcal S$.
We want a triangle whose third term has support away from $\p$.
Let $K=\Frac(B)$. For each $n$, choose elements of $H^n(Y^{\bullet})$
whose images form a basis over $K$. They define maps
$B^{r_n}[-n]\to Y^{\bullet}$ in $\D(B)$, where
$r_n=\dim_K(H^n(Y^{\bullet})\otimes_BK)$. Taking their coproduct gives
\begin{equation}\label{eq:generic-free}
 F_B^{\bullet}=\bigoplus_nB^{r_n}[-n]\longrightarrow Y^{\bullet}\longrightarrow C^{\bullet}\longrightarrow.
\end{equation}
The family defining $F_B^\bullet$ is admissible: $Y^\bullet$ has a
lower cohomological bound and each $r_n$ is finite. For every $n$,
\[
 B^{r_n}\otimes_BK\xrightarrow{\sim}H^n(Y^\bullet)\otimes_BK.
\]
The kernel of $B^{r_n}\to H^n(Y^\bullet)$ is torsion-free, as a
submodule of a free module over the domain $B$, and is zero after
tensoring with $K$. Thus this map is injective. The cohomology
sequence of \eqref{eq:generic-free} now gives
\begin{equation}\label{eq:torsion-cone}
\begin{gathered}
 0\longrightarrow B^{r_n}\longrightarrow H^n(Y^\bullet)
    \longrightarrow H^n(C^\bullet)\longrightarrow0
 \qquad(n\in\mathbb Z),\\
 C^\bullet\otimes_BK=0,
\end{gathered}
\end{equation}
because the next map $B^{r_{n+1}}\to H^{n+1}(Y^\bullet)$ is also
injective. In particular,
\[
 Z:=\Supp_R C^\bullet
    \subseteq\Supp_RY^\bullet\setminus\{\p\}
    \subseteq\Supp\mathcal S,\qquad d(Z)\le d-1.
\]
For every $\q\in Z$, induction gives $R/\q\in
\mathcal S\cap\D^b_{\fg}(R)$. Lemma~\ref{lem:prime-filtration} and
Theorem~\ref{thm:geometric} yield
\[
 C^\bullet\in\Lex(\D^b_{\fg,Z}(R))\subseteq\mathcal U.
\]
We must first show that $B\in\mathcal S$; it is not yet a generator
of $\mathcal U$. The triangle gives
\[
 C^\bullet\in\mathcal U\subseteq\mathcal S,\quad
 Y^\bullet\in\mathcal S
 \quad\Longrightarrow\quad F_B^\bullet\in\mathcal S.
\]
Since $Y^\bullet\otimes_BK\ne0$, some $r_n$ is positive, and
$B[-n]$ is a retract of $F_B^\bullet$. Therefore
\[
 B\in\mathcal S\cap\D^b_{\fg}(R),\qquad
 F_B^\bullet\in\mathcal U,\qquad Y^\bullet\in\mathcal U.
\]
Here $F_B^\bullet\in\mathcal U$ follows from the admissible coproduct
in \eqref{eq:generic-free}; the triangle then gives
$Y^\bullet\in\mathcal U$. For each $\p\in\Supp\mathcal S$, choose
$X^{\bullet}\in\mathcal S$ with $X_\p^{\bullet}\ne0$. This also proves
\eqref{eq:capture}, completing the induction.
\end{proof}

\begin{corollary}\label{cor:positive-finite}
Let $\mathcal S\subseteq\D^+_{\fg}(R)$ be a relatively localizing thick
Hom-submodule. If $d(\Supp\mathcal S)<\infty$, then
\[
 \mathcal S=\D^+_{\fg,\Supp\mathcal S}(R)
 =\Lhom(\mathcal S\cap\D^b_{\fg}(R)).
\]
\end{corollary}
\begin{proof}
Set $W=\Supp\mathcal S$. Theorem~\ref{thm:coheight} and
Lemma~\ref{lem:prime-filtration} imply
$\D^b_{\fg,W}(R)\subseteq\mathcal S\cap\D^b_{\fg}(R)$. Hence
\[
 \begin{aligned}
 \D^+_{\fg,W}(R)&=\Lex(\D^b_{\fg,W}(R))\\
 &\subseteq\Lhom(\mathcal S\cap\D^b_{\fg}(R))\\
 &\subseteq\mathcal S\subseteq\D^+_{\fg,W}(R).
 \end{aligned}
\]
The first equality is Theorem~\ref{thm:geometric}.
\end{proof}

\section{Torsion complexes over quotient domains}

\begin{lemma}\label{lem:ordinary-hom}
For every specialization-closed $W\subseteq\Spec R$,
\begin{equation}\label{eq:ordinary-hom}
 \Lhom(\D^b_{\fg,W}(R))=\Lex(\D^b_{\fg,W}(R)).
\end{equation}
\end{lemma}
\begin{proof}
Set $\mathcal L=\Lex(\D^b_{\fg,W}(R))$. It suffices to show that
$\mathcal L$ is closed under the Hom action. We first check this on
$E^{\bullet}\in\D^b_{\fg,W}(R)$. The support
$Z=\Supp E^{\bullet}$ is closed and contained in $W$. Lemma~\ref{lem:negative-model}
gives a bounded representative $Q^{\bullet}$ of finite $Z$-supported modules.
For $A^{\bullet}\in\D^-_{\fg}(R)$, choose a bounded-above representative $P^{\bullet}$ of finite
projective modules. The ordinary Hom complex $\Hom_R(P^{\bullet},Q^{\bullet})$ computes
$\RHom_R(A^{\bullet},E^{\bullet})$. It is bounded below, and each term is a finite sum
\[
 \Hom_R(P^{\bullet},Q^{\bullet})^n=\bigoplus_{i:\,Q^{i+n}\ne0}\Hom_R(P^i,Q^{i+n})
\]
of finite $Z$-supported modules. Lemma~\ref{lem:blocks} gives
$\RHom_R(A^{\bullet},E^{\bullet})\in\mathcal L$.

For fixed $A^{\bullet}$, the class
$\{X^{\bullet}\in\D^+_{\fg}(R):\RHom_R(A^{\bullet},X^{\bullet})\in\mathcal L\}$ is thick and closed under
admissible coproducts, by Lemma~\ref{lem:hom-localization}. It contains $\D^b_{\fg,W}(R)$, so it
contains $\mathcal L$. Thus $\mathcal L$ is a Hom-submodule, which proves
\eqref{eq:ordinary-hom}.
\end{proof}

\begin{proposition}\label{prop:I-X}
Let $\mathcal S\subseteq\D^+_{\fg}(R)$ be a relatively localizing thick
Hom-submodule, set $\mathcal U=\Lhom(\mathcal S\cap\D^b_{\fg}(R))$,
and fix $X^{\bullet}\in\mathcal S$. Define
\begin{equation}\label{eq:I-X}
 \mathcal I_{X^{\bullet}}=\{A^{\bullet}\in\D^-_{\fg}(R):\RHom_R(A^{\bullet},X^{\bullet})\in\mathcal U\}.
\end{equation}
There is a specialization-closed $V_{X^{\bullet}}$ with
\begin{equation}\label{eq:V-X}
 \mathcal I_{X^{\bullet}}=\D^-_{\fg,V_{X^{\bullet}}}(R),\qquad
 \{\p:\dim(R/\p)<\infty\}\subseteq V_{X^{\bullet}}.
\end{equation}
Moreover $X^{\bullet}\in\mathcal U$ if and only if $V_{X^{\bullet}}=\Spec R$.
\end{proposition}
\begin{proof}
Exactness makes $\mathcal I_{X^{\bullet}}$ thick. If
$(A_\lambda^{\bullet})_\lambda$ is an admissible family in
$\D^-_{\fg}(R)$, Lemma~\ref{lem:hom-localization} gives
\[
 \RHom_R\!\left(\bigoplus_\lambda A_\lambda^{\bullet},X^{\bullet}\right)
 \simeq\prod_\lambda\RHom_R(A_\lambda^{\bullet},X^{\bullet})
 \simeq\bigoplus_\lambda\RHom_R(A_\lambda^{\bullet},X^{\bullet}),
\]
and the family on the right is admissible. Hence
$\mathcal I_{X^{\bullet}}$ is relatively localizing. The adjunction
isomorphism
\[
 \RHom_R(A^{\bullet}\Lotimes B^{\bullet},X^{\bullet})\simeq\RHom_R(B^{\bullet},\RHom_R(A^{\bullet},X^{\bullet})),
\]
and closure of $\mathcal U$ under the Hom action show that
$\mathcal I_{X^{\bullet}}$ is a tensor ideal.
Proposition~\ref{thm:negative} gives
$\mathcal I_{X^{\bullet}}=\D^-_{\fg,V_{X^{\bullet}}}(R)$.
By Theorem~\ref{thm:coheight}, $R/\p\in\mathcal I_{X^{\bullet}}$ whenever
$\dim(R/\p)<\infty$, proving the inclusion in \eqref{eq:V-X}.
Finally,
\[
 X^{\bullet}\in\mathcal U
 \quad\Longleftrightarrow\quad
 R\in\mathcal I_{X^{\bullet}}
 \quad\Longleftrightarrow\quad
 V_{X^{\bullet}}=\Spec R.
\]
\end{proof}

\begin{condition}
\label{cond:torsion}
For every $\p\in\Spec R$, set $B=R/\p$ and $K=\Frac(B)$.
For every $C^{\bullet}\in\D^+_{\fg}(B)$ such that $C^{\bullet}\otimes_BK=0$, require
\begin{equation}\label{eq:domain-condition}
 \res C^{\bullet}\in\Loc_{\mathrm{ex},R}\bigl(\D^b_{\fg,\Supp_R(\res C^{\bullet})}(R)\bigr).
\end{equation}
Here $\res:\D(B)\to\D(R)$ is restriction of scalars along $R\to B$.
\end{condition}

\begin{theorem}\label{thm:equivalence}
For a fixed commutative noetherian ring $R$, the following are equivalent.
\begin{enumerate}[label=(\roman*),leftmargin=*]
\item Conjecture~\ref{conj:positive} holds.
\item For every specialization-closed $W\subseteq\Spec R$,
\begin{equation}\label{eq:all-geometry}
 \D^+_{\fg,W}(R)=\Lex(\D^b_{\fg,W}(R)).
\end{equation}
\item Condition~\ref{cond:torsion} holds.
\end{enumerate}
When these conditions hold, every relatively localizing thick
Hom-submodule $\mathcal S\subseteq\D^+_{\fg}(R)$ satisfies
\begin{equation}\label{eq:all-classification}
 \mathcal S=\D^+_{\fg,\Supp\mathcal S}(R).
\end{equation}
\end{theorem}
\begin{proof}
\emph{(i)$\Rightarrow$(ii).} The category $\D^+_{\fg,W}(R)$ is a relatively
localizing thick Hom-submodule, by \eqref{eq:localization} and support
preservation. Its intersection with $\D^b_{\fg}(R)$ is
$\D^b_{\fg,W}(R)$. Apply (i) and
Lemma~\ref{lem:ordinary-hom}.

\emph{(ii)$\Rightarrow$(iii).} Restriction of $C^{\bullet}$ lies in
$\D^+_{\fg,\Supp_R(\res C^{\bullet})}(R)$, so (ii) gives \eqref{eq:domain-condition}.

\emph{(iii)$\Rightarrow$(i).} Fix $\mathcal S$ and set
$\mathcal U=\Lhom(\mathcal S\cap\D^b_{\fg}(R))$. There may be no finite
bound on $\dim(R/\p)$, so the previous induction cannot be used directly.
We instead use the ascending chain condition on prime ideals. For a
prime $\p$, consider the following two assertions:
\begin{align}
 \p\in\Supp\mathcal S&\ \Longrightarrow\ R/\p\in\mathcal S\cap\D^b_{\fg}(R),
 \label{eq:good-a}\\
 \RHom_R(R/\p,X^{\bullet})&\in\mathcal U\qquad\text{for every }X^{\bullet}\in\mathcal S.
 \label{eq:good-b}
\end{align}
Suppose one of these assertions fails at some prime. Choose a
maximal such $\p$; this is possible because $R$ is noetherian.
Both assertions hold for every $\q\supsetneq\p$.
Set $B=R/\p$, $K=\Frac(B)$ and, for $X^{\bullet}\in\mathcal S$,
$Y^{\bullet}=\RHom_R(B,X^{\bullet})\in\D^+_{\fg}(B)$.

If $X_\p^{\bullet}=0$, then $Y^{\bullet}\otimes_BK=0$. All primes in $Z=\Supp_R Y^{\bullet}$ strictly
contain $\p$ and belong to $\Supp\mathcal S$. By
\eqref{eq:good-a} for these primes and Lemma~\ref{lem:prime-filtration},
$\D^b_{\fg,Z}(R)\subseteq\mathcal S\cap\D^b_{\fg}(R)$. Condition~\ref{cond:torsion}
therefore gives $Y^{\bullet}\in\Lex(\D^b_{\fg,Z}(R))\subseteq\mathcal U$.

If $X_\p^{\bullet}\ne0$, construct \eqref{eq:generic-free} over $B$, exactly as in
Theorem~\ref{thm:coheight}. The cohomology maps are injective, so its cone
$C^{\bullet}$ satisfies \eqref{eq:torsion-cone}. Moreover
\[
 Z=\Supp_R C^{\bullet}\subseteq V(\p)\setminus\{\p\}
 \subseteq\Supp\mathcal S.
\]
By \eqref{eq:good-a} at the primes strictly above $\p$, we have
$\D^b_{\fg,Z}(R)\subseteq\mathcal S\cap\D^b_{\fg}(R)$, and
Condition~\ref{cond:torsion} gives $C^{\bullet}\in\mathcal U$.
Since $Y^{\bullet}\in\mathcal S$, this implies $F_B^{\bullet}\in\mathcal S$.
By \eqref{eq:detection}, some $r_n>0$; a shift of $B$ is a retract of
$F_B^{\bullet}$, and hence $B\in\mathcal S\cap\D^b_{\fg}(R)$. Now $F_B^{\bullet}\in\mathcal U$ and
$Y^{\bullet}\in\mathcal U$.
Thus \eqref{eq:good-b} holds at $\p$. If $\p\in\Supp\mathcal S$,
choose $X^{\bullet}\in\mathcal S$ with $X_\p^{\bullet}\ne0$; the same argument
gives \eqref{eq:good-a}. This contradicts the choice of $\p$.
Both assertions therefore hold at every prime.

It remains to show that $X^\bullet\in\mathcal U$ for every
$X^\bullet\in\mathcal S$. We use the tensor ideal from
Proposition~\ref{prop:I-X}. Equation~\eqref{eq:good-b} gives
\[
 R/\p\in\mathcal I_{X^\bullet}\quad(\p\in\Spec R),\qquad
 \Supp\mathcal I_{X^\bullet}
    \supseteq\bigcup_{\p\in\Spec R}V(\p)=\Spec R.
\]
By Proposition~\ref{prop:I-X},
\[
 \mathcal I_{X^\bullet}=\D^-_{\fg}(R),\qquad
 R\in\mathcal I_{X^\bullet},\qquad
 X^\bullet=\RHom_R(R,X^\bullet)\in\mathcal U.
\]
This proves (i). Finally, \eqref{eq:good-a} and
Lemma~\ref{lem:prime-filtration} give
\[
 \D^b_{\fg,\Supp\mathcal S}(R)
       \subseteq\mathcal S\cap\D^b_{\fg}(R).
\]
Using (ii),
\[
 \D^+_{\fg,\Supp\mathcal S}(R)
   =\Lex(\D^b_{\fg,\Supp\mathcal S}(R))
   \subseteq\mathcal U\subseteq\mathcal S
   \subseteq\D^+_{\fg,\Supp\mathcal S}(R),
\]
which proves \eqref{eq:all-classification}.
\end{proof}

\begin{corollary}\label{cor:finite-torsion-test}
If $\p\in\Spec R$ satisfies $\dim(R/\p)<\infty$, then
\eqref{eq:domain-condition} holds for every
$C^{\bullet}\in\D^+_{\fg}(R/\p)$ with $C^{\bullet}\otimes_{R/\p}\Frac(R/\p)=0$.
\end{corollary}
\begin{proof}
For such a $\p$, $\Supp_R(\res C^{\bullet})\subseteq V(\p)$ has finite dimension.
Apply Theorem~\ref{thm:geometric}.
\end{proof}

\section{Duality and rings of finite character}\label{sec:infinite}

A noetherian domain $A$ has \emph{finite character} if each nonzero
element belongs to only finitely many maximal ideals. Let $K=\Frac(A)$
and let $B$ be a noetherian $A$-algebra. We call a $B$-module or a
derived object \emph{$A$-torsion} if its localization over $K$ is zero.
For an object with finite cohomology in each degree, this is equivalent
to having support in
\[
 \Sigma_A=\{\q\in\Spec B:\q\cap A\ne(0)\},
\]
where $\q\cap A$ denotes the inverse image under $A\to B$.
We write $B_\m=B\otimes_AA_\m$ and $B_K=B\otimes_AK$, and use the same
notation for complexes and supports.

\begin{lemma}\label{lem:coresolution}
Let $W\subseteq\Spec R$ be specialization closed. Suppose that, for every
integer $b$ and every
$Z^{\bullet}\in\D^+_{\fg,W}(R)\cap\D^{\geq b}(R)$, there are a finite $W$-supported module
$N$ and a map $Z^{\bullet}\to N[-b]$ injective on $H^b$.
For every integer $a$, each $X^{\bullet}\in\D^+_{\fg,W}(R)\cap\D^{\geq a}(R)$
then has a representative $Q^{\bullet}$ with
finite $W$-supported terms and $Q^i=0$ for $i<a$. In particular,
$\D^+_{\fg,W}(R)=\Lex(\D^b_{\fg,W}(R))$.
\end{lemma}
\begin{proof}
After a shift, assume $a=0$. We construct the terms of the representative
one at a time. Set $X_0^{\bullet}=X^{\bullet}$ and choose triangles
\begin{equation}\label{eq:coresolution-fibers}
 X_{n+1}^{\bullet}\longrightarrow X_n^{\bullet}\xrightarrow{f_n}N_n[-n]
 \longrightarrow X_{n+1}^{\bullet}[1],\qquad X_n^{\bullet}\in\D^{\geq n}(R).
\end{equation}
Since $H^n(f_n)$ is injective, $X_{n+1}^{\bullet}\in\D^{\geq n+1}(R)$.
We will arrange these choices into compatible triangles
\begin{equation}\label{eq:finite-heads}
 X_{n+1}^{\bullet}\longrightarrow X^{\bullet}\xrightarrow{\alpha_n}P_n^{\bullet}\longrightarrow X_{n+1}^{\bullet}[1],
\end{equation}
where $P_n^{\bullet}$ is the complex $N_0\to\cdots\to N_n$ in degrees
$0,\ldots,n$. For $n=0$, use \eqref{eq:coresolution-fibers}. Suppose that $P_n^{\bullet}$
has been constructed. Its next differential comes from the composite
\[
 P_n^{\bullet}\longrightarrow X_{n+1}^{\bullet}[1]\longrightarrow N_{n+1}[-n].
\]
Since $P_n^{\bullet}\in\D^{\leq n}(R)$,
\[
 \Hom_{\D(R)}(P_n^{\bullet},N_{n+1}[-n])
 \cong\Hom_R(H^n(P_n^{\bullet}),N_{n+1}).
\]
The composite is therefore represented by a module map
$d_n:N_n\to N_{n+1}$ which factors through $\coker d_{n-1}$
(where $d_{-1}=0$). In particular $d_nd_{n-1}=0$. Its fiber is
\[
 P_{n+1}^{\bullet}=(N_0\longrightarrow\cdots\longrightarrow N_n
 \xrightarrow{d_n}N_{n+1}).
\]
Apply the octahedral axiom to
$X_{n+2}^{\bullet}\to X_{n+1}^{\bullet}\to X^{\bullet}$.
It gives \eqref{eq:finite-heads} for $n+1$, with
$\alpha_n=\pi_n\alpha_{n+1}$, where
$\pi_n:P_{n+1}^{\bullet}\to P_n^{\bullet}$ is the brutal-truncation projection.

Set $Q^{\bullet}=(N_0\xrightarrow{d_0}N_1\xrightarrow{d_1}\cdots)$.
It remains to identify $Q^{\bullet}$ with $X^{\bullet}$.
The maps $\pi_n$ are degreewise surjective and eventually constant
in each degree. Thus, in $\D(\Mod R)$,
\[
 \operatorname{holim}_nP_n^{\bullet}\simeq\varprojlim_nP_n^{\bullet}=Q^{\bullet}.
\]
Here the homotopy inverse limit is defined by the ambient triangle
\[
 \operatorname{holim}_nP_n^{\bullet}\longrightarrow\prod_nP_n^{\bullet}
 \xrightarrow{1-\mathrm{shift}}\prod_nP_n^{\bullet}\longrightarrow
 (\operatorname{holim}_nP_n^{\bullet})[1],
\]
where $\mathrm{shift}$ is induced by the transition maps $\pi_n$.
The compatible maps $\alpha_n$ define a map
$X^{\bullet}\to\prod_nP_n^{\bullet}$ whose composite with
$1-\mathrm{shift}$ is zero. The exact Hom sequence of this triangle
therefore gives a map $X^{\bullet}\to\operatorname{holim}_nP_n^{\bullet}$.
For a fixed $i$ and $n\geq i+1$, triangle
\eqref{eq:finite-heads} gives
\[
 H^i(\alpha_n):H^i(X^\bullet)\xrightarrow{\sim}H^i(P_n^\bullet).
\]
The inverse systems of cohomology are eventually constant. The Milnor
sequence for the displayed product triangle is
\[
 0\longrightarrow\varprojlim\nolimits_n^1 H^{i-1}(P_n^\bullet)
 \longrightarrow H^i(\operatorname{holim}_nP_n^\bullet)
 \longrightarrow\varprojlim_nH^i(P_n^\bullet)\longrightarrow0.
\]
The left term is zero, and the right term is $H^i(X^\bullet)$.
Hence $X^\bullet\to Q^\bullet$ is a quasi-isomorphism.
Each $Q^i=N_i$ is finite and supported in $W$.
Lemma~\ref{lem:blocks} now gives the asserted equality of subcategories.
\end{proof}

\begin{lemma}
\label{lem:finite-id-localization}
Let $X^{\bullet}\in\D^+_{\fg}(R)$ and let $E^{\bullet}\in\D^b_{\fg}(R)$ admit a bounded
injective representative. For every multiplicative subset $\Phi\subseteq R$,
\begin{equation}\label{eq:finite-id-localization}
 \Hom_{\D(R)}(X^{\bullet},E^{\bullet})_\Phi
 \cong\Hom_{\D(R_\Phi)}(X_\Phi^{\bullet},E_\Phi^{\bullet}).
\end{equation}
\end{lemma}
\begin{proof}
A bounded injective representative of the target lets us discard
all sufficiently high degrees of the source. Choose such a representative
$I^{\bullet}$, with $I^j=0$ for $j>v$. For $n>v+1$,
\[
 \Hom_{\D(R)}(\tau_{\geq n+1}X^{\bullet},E^{\bullet}[u])=0\qquad(u=0,1),
\]
because these groups are computed by maps to $I^{\bullet}[u]$. Thus
\[
 \Hom_{\D(R)}(X^{\bullet},E^{\bullet})
 \cong\Hom_{\D(R)}(\tau_{\leq n}X^{\bullet},E^{\bullet}).
\]
Localization preserves injectives over a noetherian ring
\cite[Tag~0A6I]{Stacks}, so the same truncation works over $R_\Phi$.
The source $\tau_{\leq n}X^{\bullet}$ is bounded, so
\eqref{eq:bounded-hom-localization} gives the result.
\end{proof}

\begin{lemma}\label{lem:koszul-embedding}
Let $E^{\bullet}\in\D^b_{\fg}(R)$, let $L\subseteq H^a(E^{\bullet})$, and let
$J=(f_1,\ldots,f_c)$ annihilate $L$. There are positive integers
$n_1,\ldots,n_c$ such that the canonical map
\[
 E^{\bullet}\longrightarrow E^{\bullet}\otimes_R^{\mathbf L}
 K_R^{\bullet}(f_1^{n_1},\ldots,f_c^{n_c})
\]
is injective on $L$ in degree $a$. The target is supported in
$\Supp E^{\bullet}\cap V(J)$. If $E^{\bullet}$ has finite injective dimension, so does
the target. Here
\[
 K_R^{\bullet}(f_1,\ldots,f_c)=
 \bigotimes_{j=1}^c(R\xrightarrow{f_j}R),
\]
with each factor in degrees $-1,0$ and with canonical map
$R\to K_R^{\bullet}(f_j)$ given by inclusion in degree zero.
\end{lemma}
\begin{proof}
First take one element $f$. The triangle
$E^{\bullet}\xrightarrow{f^n}E^{\bullet}\to E^{\bullet}\otimes_RK_R^{\bullet}(f^n)\to E^{\bullet}[1]$ gives
\[
 \ker\bigl(H^a(E^{\bullet})\to H^a(E^{\bullet}\otimes_RK_R^{\bullet}(f^n))\bigr)=f^nH^a(E^{\bullet}).
\]
Artin--Rees gives an integer $d$ such that, for $n>d$,
\[
 L\cap f^nH^a(E^{\bullet})
 =f^{n-d}\bigl(L\cap f^dH^a(E^{\bullet})\bigr)=0,
\]
since $fL=0$ \cite[Tag~00IN]{Stacks}. Thus the map to the Koszul tensor
is injective on $L$. Repeat this argument for $f_1,\ldots,f_c$.
At each step the image of $L$ is still annihilated by $J$.
Outside $V(J)$ one Koszul factor is contractible, which gives the
support assertion. Finally, tensoring a bounded injective complex
with a bounded finite-free complex again gives a bounded complex
of injectives.
\end{proof}

A dualizing complex is an object $\omega^{\bullet}\in\D^b_{\fg}(R)$ of finite
injective dimension for which
$R\to\RHom_R(\omega^{\bullet},\omega^{\bullet})$ is an isomorphism. The functor
$D_{\omega^{\bullet}}=\RHom_R(-,\omega^{\bullet})$ exchanges $\D^-_{\fg}(R)$ and
$\D^+_{\fg}(R)$, with $D_{\omega^{\bullet}}^2\simeq1$
\cite[Tag~0A7C]{Stacks}.

\begin{lemma}
\label{lem:dualizing-supported-model}
Suppose that $R$ admits a dualizing complex. For every specialization-closed
$W\subseteq\Spec R$ and every
$X^{\bullet}\in\D^+_{\fg,W}(R)\cap\D^{\geq a}(R)$, there is a termwise finite
$W$-supported representative of $X^{\bullet}$, zero below $a$.
\end{lemma}
\begin{proof}
By Lemma~\ref{lem:coresolution}, it suffices to embed $H^a(X^{\bullet})$
in a finite $W$-supported module by a map defined on $X^{\bullet}$.
If $H^a(X^{\bullet})=0$, use the zero module. Otherwise choose a
bounded-above finite-free representative $P^{\bullet}$ of
$D_{\omega^{\bullet}}X^{\bullet}$. Choose $n$ such that
$n-1\leq\inf\omega^{\bullet}-a-2$, and let $(P^{\bullet})^{\geq n}$ be its brutal lower
truncation. The cone $T^{\bullet}$ of $(P^{\bullet})^{\geq n}\to D_{\omega^{\bullet}} X^{\bullet}$ satisfies
$\sup T^{\bullet}\leq n-1$. Since
\[
 \inf D_{\omega^{\bullet}} T^{\bullet}\geq\inf\omega^{\bullet}-\sup T^{\bullet}\geq a+2,
\]
duality gives a map
\begin{equation}\label{eq:dualizing-head}
 X^{\bullet}\longrightarrow E^{\bullet}=D_{\omega^{\bullet}}((P^{\bullet})^{\geq n})
\end{equation}
which is an isomorphism on $H^a$.
The object $E^{\bullet}$ is bounded, has finite cohomology, and has
finite injective dimension. Its support may be too large.
Set $J=\operatorname{ann}_R H^a(X^{\bullet})$ and apply
Lemma~\ref{lem:koszul-embedding} to \eqref{eq:dualizing-head}. This gives a map
$X^{\bullet}\to E^{\prime\bullet}$ injective on $H^a$, with
$E^{\prime\bullet}\in\D^b_{\fg,V(J)}(R)$.

By Lemma~\ref{lem:negative-model} and good truncation,
$\tau_{\geq a}E^{\prime\bullet}$ has a bounded representative $Q^{\bullet}$ of finite
$V(J)$-supported modules, zero below $a$. The composite
\[
 X^{\bullet}\longrightarrow E^{\prime\bullet}\longrightarrow\tau_{\geq a}E^{\prime\bullet}
 \simeq Q^{\bullet}\longrightarrow Q^a[-a]
\]
is injective on $H^a$. The finite module $Q^a$ has support in
$V(J)=\Supp_RH^a(X^{\bullet})\subseteq W$.
Lemma~\ref{lem:coresolution} applies.
\end{proof}

\begin{theorem}
\label{prop:dualizing-gluing}
Let $A$ be a noetherian domain of finite character, and let $B$ be a
noetherian $A$-algebra. Suppose there is $\omega^{\bullet}\in\D^b_{\fg}(B)$ such
that $\omega_\m^{\bullet}$ is dualizing over $B_\m$ for every
$\m\in\Max A$ with $B_\m\ne0$. Then, for every specialization-closed
$W\subseteq\Spec B$,
\begin{equation}\label{eq:local-duality-generation}
 \D^+_{\fg,W}(B)=\Lex(\D^b_{\fg,W}(B)).
\end{equation}
If $W\subseteq\Sigma_A$, every object of $\D^+_{\fg,W}(B)$ admits a
bounded-below representative by finite $W$-supported modules.
\end{theorem}
\begin{proof}
Set $K=\Frac(A)$. The injective dimensions of the local dualizing
complexes need not have a common bound. We only need such a bound
after tensoring with an $A$-torsion perfect complex. In fact,
\begin{equation}\label{eq:torsion-dualizing-target}
 \begin{gathered}
 Q^{\bullet}\in\Perf(B),\qquad Q_K^{\bullet}=0\\
 \Longrightarrow\quad
 \omega^{\bullet}\otimes_B^{\mathbf L}Q^{\bullet}
 \text{ has finite injective dimension over }B.
 \end{gathered}
\end{equation}
Set $E^{\bullet}=\omega^{\bullet}\otimes_B^{\mathbf L}Q^{\bullet}$.
Since $Q^{\bullet}$ is perfect, $E^{\bullet}\in\D^b_{\fg}(B)$.
Choose $0\ne s\in A$ annihilating all the finitely many nonzero
cohomology modules of $Q^{\bullet}$, and write
\[
 \{\m\in\Max A:s\in\m\}=\{\m_1,\ldots,\m_r\}.
\]
If $r=0$, then $Q^{\bullet}=0$. Otherwise each
$E^{\bullet}_{\m_j}$ has finite injective dimension, since it is
the tensor product of a dualizing complex with a perfect complex.
There are only finitely many $\m_j$, so we can choose integers
$a,b$ such that
\[
 E^{\bullet}\in\D^{\geq a}(B),\qquad
 \Ext^i_{B_{\m_j}}(M,E^{\bullet}_{\m_j})=0
 \quad(i>b, M\in\Mod B_{\m_j},\ 1\leq j\leq r).
\]
At every other maximal ideal $\m$ of $A$, both
$Q^{\bullet}_{\m}$ and $E^{\bullet}_{\m}$ vanish. To pass back
to $B$, take any ideal $I\subseteq B$. Since $B/I$ is finite and
concentrated in degree zero,
\[
 \Ext^i_B(B/I,E^{\bullet})_{\m}
 \cong
 \Ext^i_{B_{\m}}(B_{\m}/IB_{\m},E^{\bullet}_{\m}).
\]
Thus $\Ext^i_B(B/I,E^{\bullet})_{\m}=0$ for all
$\m\in\Max A$ and $i>b$, so
$\Ext^i_B(B/I,E^{\bullet})=0$ for $i>b$. The same group vanishes
for $i<a$ by the standard $t$-structure. The injective-amplitude
criterion \cite[Tag~0A5T]{Stacks} proves
\eqref{eq:torsion-dualizing-target}.

We next treat $W\subseteq\Sigma_A$. By Lemma~\ref{lem:coresolution},
we have to embed the lowest cohomology of
$Z^{\bullet}\in\D^+_{\fg,W}(B)\cap\D^{\geq a}(B)$ in a finite
$W$-supported module. If $H^a(Z^{\bullet})=0$, take the zero module.
Otherwise $H^a(Z^{\bullet})$ is finite over $B$ and $A$-torsion;
a nonzero element of $A$ kills all of it. Thus only finitely many
$\m\in\Max A$ satisfy $H^a(Z^{\bullet})_\m\ne0$.
Fix one of them. The construction
\eqref{eq:dualizing-head} over $B_\m$ gives
\[
 Z_\m^{\bullet}\longrightarrow
 \RHom_{B_\m}(F_\m^{\bullet},\omega_\m^{\bullet})
 \quad\text{injective on }H^a,
\]
where $F_\m^{\bullet}$ is a bounded complex of finite free $B_\m$-modules.

We need a bounded finite-free $B$-complex whose localization is
$F_\m^{\bullet}$. Choose $u\in A\setminus\m$ clearing the
finitely many differential matrices, and lift
$u d^i_{F_\m^{\bullet}}$ to matrices $D_i$ over $B$. Their adjacent
products are zero after localization. Choose $v\in A\setminus\m$
such that $vD_{i+1}D_i=0$ for every $i$. The matrices $vD_i$
now define a complex $F^{\bullet}$. Its localized differentials
are $uvd^i_{F_\m^\bullet}$, and
\[
 \psi^i=(uv)^{-i}\operatorname{id},\qquad
 \psi^{i+1}(uvd^i_{F_\m^\bullet})
       =d^i_{F_\m^\bullet}\psi^i
\]
gives an isomorphism $F_\m^\bullet\simeq(F^\bullet)_\m$.

Set $E^{0,\bullet}=\RHom_B(F^{\bullet},\omega^{\bullet})$ and
$J=\operatorname{ann}_B H^a(Z^{\bullet})$. Then
\[
 V(J)\subseteq W,\qquad J\cap A\ne(0).
\]
Choose generators $f_1,\ldots,f_c$ of $J$.
Lemma~\ref{lem:koszul-embedding}, applied over $B_\m$, gives
positive integers $n_j$ such that the local composite into
\[
 E^{\prime\bullet}=E^{0,\bullet}\otimes_BK_B^{\bullet}(f_1^{n_1},\ldots,f_c^{n_c})
\]
remains injective on $H^a$.
The perfect complex
\[
 (F^{\bullet})^\vee\otimes_BK_B^{\bullet}(f_1^{n_1},\ldots,f_c^{n_c}),
 \qquad (F^{\bullet})^\vee=\Hom_B(F^{\bullet},B),
\]
is supported in $V(J)$ and becomes zero over $K$.
By \eqref{eq:torsion-dualizing-target}, $E^{\prime\bullet}$ has
finite injective dimension over $B$. We can therefore apply
Lemma~\ref{lem:finite-id-localization}: after multiplication by
an element of $A\setminus\m$, the local map comes from a map
$Z^{\bullet}\to E^{\prime\bullet}$. Multiplication by this
element does not affect injectivity after localization at $\m$.

Now choose a bounded representative of $\tau_{\geq a}E^{\prime\bullet}$ by finite
$V(J)$-supported modules, zero below $a$ (Lemma~\ref{lem:negative-model}).
Projection to its lowest term gives
\[
 Z^{\bullet}\longrightarrow N(\m)[-a]
\]
injective on $H^a$ after localization at $\m$, with $N(\m)$ finite
and $W$-supported. Take the sum of these maps over the finitely many maximal ideals
where $H^a(Z^{\bullet})$ is nonzero. The kernel on $H^a$
vanishes at every maximal ideal of $A$, hence is zero.
Lemma~\ref{lem:coresolution} gives the required representative and
\begin{equation}\label{eq:A-torsion-models}
 \D^+_{\fg,W}(B)=\Lex(\D^b_{\fg,W}(B))
 \qquad(W\subseteq\Sigma_A).
\end{equation}

It remains to allow arbitrary $W$. Take $X^{\bullet}\in\D^+_{\fg,W}(B)$.
If $B_K=0$, then $W\subseteq\Sigma_A$ and \eqref{eq:A-torsion-models}
applies. If $B_K\ne0$, then $B_\m\ne0$ for every $\m\in\Max A$.
Choose any such $\m$ and localize the dualizing
complex $\omega_\m^{\bullet}$ further to $B_K$. Then $\omega_K^{\bullet}$ is dualizing
\cite[Tag~0A7G]{Stacks}. Lemma~\ref{lem:dualizing-supported-model}
gives a termwise finite $W_K$-supported positive representative of $X_K^{\bullet}$.
By Lemma~\ref{lem:scalar-spreading}, there is $Y^{\bullet}\to X^{\bullet}$ with
$Y^{\bullet}\in\Lex(\D^b_{\fg,W}(B))$ which becomes an isomorphism over $K$.
The cone is supported in $W\cap\Sigma_A$, so
\eqref{eq:A-torsion-models} applies to it. The resulting triangle
shows that $X^{\bullet}\in\Lex(\D^b_{\fg,W}(B))$.
\end{proof}

\begin{corollary}
\label{thm:finite-character}
Let $A$ be a noetherian domain of finite character and let $B$ be a
commutative $A$-algebra of finite type. Assume that either $A$ or $B$
is regular. For every specialization-closed $W\subseteq\Spec B$,
\begin{equation}\label{eq:finite-character-geometry}
 \D^+_{\fg,W}(B)=\Lex(\D^b_{\fg,W}(B)).
\end{equation}
Every relatively localizing thick Hom-submodule
$\mathcal S\subseteq\D^+_{\fg}(B)$ satisfies
\begin{equation}\label{eq:finite-character-telescope}
 \mathcal S=\D^+_{\fg,\Supp\mathcal S}(B)
 =\Lhom(\mathcal S\cap\D^b_{\fg}(B)).
\end{equation}
\end{corollary}
\begin{proof}
Write $B$ as a quotient of $P=A[t_1,\ldots,t_r]$.
For $\m\in\Max A$,
\[
 \dim B_\m\leq\dim P_\m=\dim A_\m+r<\infty
 \quad\text{\cite[Tag~02IJ]{Stacks}}.
\]
If $B$ is regular, each $B_\m$ is regular of finite dimension,
hence has finite global dimension \cite[Tag~00OE]{Stacks}.
Apply Theorem~\ref{prop:dualizing-gluing} with $\omega^{\bullet}=B$.

If $A$ is regular, then $P$ is regular. Every finite module over a
regular noetherian ring is perfect \cite[Tag~066Z]{Stacks}, so
\[
 \omega^{\bullet}=\RHom_P(B,P)\in\D^b_{\fg}(B).
\]
Each $P_\m$ is regular of finite dimension, hence is a dualizing complex
over itself. Localization and finite coinduction show that
\[
 \omega_\m^{\bullet}\simeq\RHom_{P_\m}(B_\m,P_\m)
\]
is dualizing over each nonzero $B_\m$ \cite[Tag~0AX0]{Stacks}.
Thus Theorem~\ref{prop:dualizing-gluing} gives
\eqref{eq:finite-character-geometry} also in this case.
Apply Theorem~\ref{thm:equivalence} over $B$ to obtain
\eqref{eq:finite-character-telescope}.
\end{proof}

\begin{example}[Nagata's infinite-dimensional domain]\label{ex:nagata}
Let $P=k[x_1,x_2,\ldots]$ and
\[
 I_i=\{2^{i-1},\ldots,2^i-1\},\quad
 \p_i=(x_j:j\in I_i),\quad
 S=P\setminus\bigcup_{i\geq1}\p_i,\quad A=S^{-1}P.
\]
Nagata's construction makes $A$ noetherian with maximal ideals
$\m_i=\p_iA$ \cite[Section~110.16, Tag~02JC]{Stacks}.
The local ring $A_{\m_i}=P_{\p_i}$ is a polynomial regular local ring over
the fraction field in the variables outside $I_i$, of dimension
$|I_i|=2^{i-1}$. Thus $A$ is regular and $\dim A=\infty$.
A nonzero polynomial involves variables in only finitely many blocks,
so it lies in only finitely many $\p_i$. The same holds for a nonzero
element of $A$ and the ideals $\m_i$. Hence $A$ has finite character.

Corollary~\ref{thm:finite-character} applies to every finite-type
$A$-algebra. For $B=A[t_1,\ldots,t_r]$, $r\geq1$, the prime
$\p=(t_1,\ldots,t_r)$ satisfies
\[
 \dim(B/\p)=\dim A=\infty.
\]
The singular algebra $B=A[t,\varepsilon]/(\varepsilon^2)$ is also
covered, with $\p=(t,\varepsilon)$ and
$\dim(B/\p)=\infty$.
\end{example}

\section{Consequences and a counterexample}

\begin{corollary}\label{prop:coinduced}
Let $A^{\bullet}\in\D^-_{\fg}(R)$, $E^{\bullet}\in\D^b_{\fg}(R)$, and $W=\Supp A^{\bullet}$. Then
\begin{equation}\label{eq:coinduced-generation}
 \Lhom(\RHom_R(A^{\bullet},E^{\bullet}))=
 \Lhom\{\RHom_R(P^{\bullet},E^{\bullet}):P^{\bullet}\in\Perf_W(R)\}.
\end{equation}
Each $\RHom_R(P^{\bullet},E^{\bullet})$ on the right belongs to $\D^b_{\fg}(R)$.
\end{corollary}
\begin{proof}
Let $\mathcal L$ be a relatively localizing thick Hom-submodule. The class
\[
 \{B^{\bullet}\in\D^-_{\fg}(R):\RHom_R(B^{\bullet},E^{\bullet})\in\mathcal L\}
\]
is a relatively localizing thick tensor ideal, by the tensor--Hom
adjunction and Lemma~\ref{lem:hom-localization}.
Proposition~\ref{thm:negative} gives
$A^{\bullet}\in\Lten(\Perf_W(R))$ and
$\Perf_W(R)\subseteq\Lten(A^{\bullet})$. Taking $\mathcal L$ first to be the left side of
\eqref{eq:coinduced-generation} and then the right side gives
the two inclusions.
For perfect $P^{\bullet}$, set $(P^{\bullet})^\vee=\RHom_R(P^{\bullet},R)$. Then
$\RHom_R(P^{\bullet},E^{\bullet})\simeq (P^{\bullet})^\vee\Lotimes E^{\bullet}\in\D^b_{\fg}(R)$.
\end{proof}

\begin{corollary}\label{prop:pieces}
Suppose a specialization-closed set has a disjoint decomposition
$W=\coprod_{\lambda\in\Lambda}W_\lambda$ into specialization-closed subsets
with $d(W_\lambda)<\infty$ for each $\lambda$. Then
\[
 \D^+_{\fg,W}(R)=\Lex(\D^b_{\fg,W}(R)).
\]
Every relatively localizing thick Hom-submodule
$\mathcal S\subseteq\D^+_{\fg}(R)$ with $\Supp\mathcal S\subseteq W$ satisfies
\[
 \mathcal S=\D^+_{\fg,\Supp\mathcal S}(R)
 =\Lhom(\mathcal S\cap\D^b_{\fg}(R)).
\]
\end{corollary}
\begin{proof}
First let $M\in\Mod_W R$ be finite. Each irreducible component
of $\Supp M$ lies in the unique $W_\lambda$ containing its generic
point. Grouping the components in this way gives a finite disjoint
closed decomposition of $\Supp M$. The corresponding idempotents
of $R/\operatorname{ann}_R M$ decompose $M$.
Every homomorphism preserves these components: an image has support
contained in the support of its source, and a submodule of a different
component has disjoint support. The decompositions are therefore
compatible with inclusions of finite submodules. Every module is the filtered union of its finite submodules. Thus,
for arbitrary $M\in\Mod_W R$,
\[
 \begin{aligned}
 M&=\bigoplus_\lambda\Gamma_{W_\lambda}M,\\
 \Gamma_{W_\lambda}M&=
 \{m\in M:\Supp_R(Rm)\subseteq W_\lambda\},
 \end{aligned}
\]
and these functors are exact on $\Mod_W R$.

By the bounded-below torsion comparison
\cite[Lemma~51.5.4, Tag~0EF2]{Stacks}, an object
$X^{\bullet}\in\D^+_{\fg,W}(R)$ has a bounded-below representative
in $\Mod_W R$. Decomposing each term gives
$X^{\bullet}\simeq\bigoplus_\lambda X_\lambda^{\bullet}$, where
$H^n(X_\lambda^{\bullet})=\Gamma_{W_\lambda}H^n(X^{\bullet})$. Each cohomology module is a
direct summand of the finite module $H^n(X^{\bullet})$, and only finitely many
$\lambda$ occur in a fixed degree. Thus the decomposition is admissible
in $\D^+_{\fg}(R)$. Theorem~\ref{thm:geometric} places every $X_\lambda^{\bullet}$ in
$\Lex(\D^b_{\fg,W_\lambda}(R))$, proving the first assertion.

If $\Supp\mathcal S\subseteq W$, each prime in $\Supp\mathcal S$ has finite
coheight: $V(\p)$ is contained in its $W_\lambda$. Theorem~\ref{thm:coheight} therefore gives
$R/\p\in\mathcal S\cap\D^b_{\fg}(R)$ for every
$\p\in\Supp\mathcal S$. Apply the first assertion to the
disjoint subsets $\Supp\mathcal S\cap W_\lambda$ and then use
Lemma~\ref{lem:prime-filtration}.
\end{proof}

\Needspace{10\baselineskip}
\begin{theorem}
\label{thm:no-finite-model}
Let $A$ be the domain in Example~\ref{ex:nagata} and set $B=A[t]$.
There is an object $C^{\bullet}\in\D^+_{\fg}(B)$ such that
\begin{enumerate}[label=(\roman*),leftmargin=*]
\item $C^{\bullet}$ has a bounded-below representative on which $t$ acts termwise
as zero, and $\Supp_BC^{\bullet}=V(t)$;
\item $C^{\bullet}\in\Lex(\D^b_{\fg,V(t)}(B))$;
\item $C^{\bullet}$ has no bounded-below representative by finite $B$-modules.
\end{enumerate}
\end{theorem}
\begin{proof}
We choose nonzero extension classes of increasing degrees at different
maximal ideals. A representative by finite modules would give a monomorphism
$A\to N$, with $N$ finite over $B$, sending all these classes to zero.
We will show that such a map cannot exist.
Take $A$, $\m_i$, and $d_i=2^{i-1}$ from Example~\ref{ex:nagata}.
Use only $i\geq3$, and set $k_i=A/\m_i$ and $B=A[t]$.
Koszul resolutions over the regular local rings $A_{\m_i}$ give
\[
 \Ext_A^{d_i}(k_i,A)_{\m_i}
 \cong\Ext_{A_{\m_i}}^{d_i}(k_i,A_{\m_i})\cong k_i.
\]
Choose $\eta_i\in\Ext_A^{d_i}(k_i,A)$ with nonzero localization.
The chosen classes define the components of a morphism
\[
 T^{\bullet}=\bigoplus_{i\geq3}k_i[-(d_i-1)]\longrightarrow A[1]
\]
in $\D(A)$. The degrees $d_i-1$ tend to infinity, so the coproduct is admissible.
Complete this morphism to a triangle
$A\to C_A^{\bullet}\to T^{\bullet}\to A[1]$.
Let $i_*:\D(A)\to\D(B)$ be restriction along $B\to B/(t)=A$, and set
$C^{\bullet}=i_*C_A^{\bullet}$. Thus
\begin{equation}\label{eq:star-extension}
 A\longrightarrow C^{\bullet}\longrightarrow i_*T^{\bullet}\xrightarrow{e}A[1],
 \qquad e_i=i_*(\eta_i).
\end{equation}
One has
\[
 H^0(C^{\bullet})=A,\qquad H^{d_i-1}(C^{\bullet})=k_i,
\]
and the other cohomology modules are zero. This gives (i).
The objects $A$ and $k_i$ are finite and supported in $V(t)$;
hence \eqref{eq:star-extension} also gives (ii).

For (iii), suppose that $C^{\bullet}\simeq Q^{\bullet}$ with
$Q^{\bullet}$ bounded below and termwise finite over $B$.
Replace $Q^{\bullet}$ by its good truncation so that $Q^j=0$
for $j<0$. The projection $Q^{\bullet}\to Q^0[0]$ restricts on
$H^0(Q^{\bullet})=A$ to the inclusion of cycles. Write this
monomorphism as $u:A\hookrightarrow N$, where $N=Q^0$ is finite
over $B$. If $\iota:A\to C^\bullet$ is the first map of
\eqref{eq:star-extension}, the projection gives
\[
 \widetilde u:C^\bullet\longrightarrow N,\qquad
 \widetilde u\iota=u,\qquad u[1]e=0.
\]
Thus exactness of Hom gives
\begin{equation}\label{eq:simultaneous-effacement}
 u_*(e_i)=0\quad\text{in }\Ext_B^{d_i}(k_i,N)\quad(i\geq3).
\end{equation}
The forgetful functor $U:\D(B)\to\D(A)$ satisfies $Ui_*=1$.
Applying it to \eqref{eq:simultaneous-effacement} gives
$u_*(\eta_i)=0$ in $\Ext_A^{d_i}(k_i,N)$.

Here $N$ is finite over $B$, but need not be finite over $A$.
We can nevertheless split $u$ after inverting one nonzero element
of $A$. Generic freeness \cite[Lemma~10.118.1, Tag~051R]{Stacks} gives
$0\ne s_0\in A$ and an $A_{s_0}$-basis $(v_\lambda)_{\lambda\in\Lambda}$
of $N_{s_0}$. Write
\[
 u_{s_0}(1)=\sum_{\lambda\in F}c_\lambda v_\lambda,
 \qquad F\subseteq\Lambda\text{ finite}.
\]
The injection $u_{s_0}:A_{s_0}\hookrightarrow N_{s_0}$ is nonzero.
Choose $\lambda_0\in F$ with $c_{\lambda_0}\ne0$ and write
$c_{\lambda_0}=c/s_0^h$, where $0\ne c\in A$. Set $s=s_0c$.
The coordinate functional
\[
 \ell_s:N_s\longrightarrow A_s,\qquad
 \ell_s(v_{\lambda_0})=c_{\lambda_0}^{-1},\quad
 \ell_s(v_\lambda)=0\quad(\lambda\ne\lambda_0),
\]
is defined on the possibly infinite-rank free module and satisfies
$\ell_su_s=1_{A_s}$. Since $s$ lies in only finitely many maximal
ideals, choose $i\geq3$ with $s\notin\m_i$. The splitting then
localizes to
\[
 A_{\m_i}\xrightarrow{u_{\m_i}}N_{\m_i}
 \xrightarrow{\ell_{\m_i}}A_{\m_i},
 \qquad\ell_{\m_i}u_{\m_i}=1.
\]
Consequently
\[
 (\ell_{\m_i})_*(u_{\m_i})_*
 =1_{\Ext_{A_{\m_i}}^{d_i}(k_i,A_{\m_i})}.
\]
To compare this with the global extension, note that $k_i$ is
perfect over the regular ring $A$. Its bounded finite-projective
resolution gives
\[
 \Ext_A^{d_i}(k_i,N)_{\m_i}
 \cong\Ext_{A_{\m_i}}^{d_i}(k_i,N_{\m_i}),
\]
although $N$ need not be finite over $A$. Hence
$u_*(\eta_i)=0$ would imply
$(u_{\m_i})_*((\eta_i)_{\m_i})=0$, and applying
$(\ell_{\m_i})_*$ would give $(\eta_i)_{\m_i}=0$.
This contradicts the choice of $\eta_i$ and proves (iii).
\end{proof}

\begin{corollary}\label{cor:comparison-failure}
For $B$ in Theorem~\ref{thm:no-finite-model}, the natural functors
\[
 \begin{aligned}
 F^+&:\D^+(B\text{-}\mathrm{mod})\longrightarrow\D^+_{\fg}(B),\\
 F&:\D(B\text{-}\mathrm{mod})\longrightarrow\D_{\fg}(B)
 \end{aligned}
\]
are neither full nor essentially surjective.
\end{corollary}
\begin{proof}
Let $C^{\bullet}$ be the object constructed in Theorem~\ref{thm:no-finite-model}.
If a complex of finite $B$-modules represented $C^{\bullet}$, its good lower
truncation at zero would be a bounded-below such representative,
since $H^{<0}(C^{\bullet})=0$. Thus $C^{\bullet}$ belongs to neither essential image.

The two ends $i_*T^{\bullet}$ and $A[1]$ of
\eqref{eq:star-extension} do have bounded-below representatives
by finite modules. If the connecting morphism
$e:i_*T^{\bullet}\to A[1]$ lifted to a morphism $\widetilde e$ in
either source category, exactness would give
\[
 F^+\bigl(\cone(\widetilde e)[-1]\bigr)\simeq C^\bullet
 \quad\text{or}\quad
 F\bigl(\cone(\widetilde e)[-1]\bigr)\simeq C^\bullet.
\]
Either isomorphism would give a representative by finite modules of
$C^\bullet$, a contradiction. Thus neither functor is full.
\end{proof}

\begin{example}\label{ex:localization-failure}
Let $k$ be a field,
$R=k[t,\varepsilon]/(\varepsilon^2)$, $M=R/(\varepsilon)$, and
$X^{\bullet}=\bigoplus_{n\geq0}M[-n]\in\D^+_{\fg}(R)$. The periodic resolution with
differential $\varepsilon$ gives $\Ext_R^n(M,M)=k[t]$ for $n\geq0$.
Thus
\[
 \begin{aligned}
 \Hom_{\D(R)}(X^{\bullet},M)&\cong\prod_{n\geq0}k[t],\\
 \Hom_{\D(R_t)}(X_t^{\bullet},M_t)&\cong\prod_{n\geq0}k[t,t^{-1}].
 \end{aligned}
\]
Under these identifications, localization is the nonsurjective map
\[
 \left(\prod_{n\geq0}k[t]\right)_t
 \longrightarrow\prod_{n\geq0}k[t,t^{-1}];
\]
the element $(t^{-n})_{n\geq0}$ has no common denominator.
\end{example}

\begin{proposition}\label{prop:tail-obstruction}
Fix a specialization-closed $W$, set $\mathcal U_W=\Lex(\D^b_{\fg,W}(R))$, and
let $q:\D^+_{\fg,W}(R)\to\mathcal Q=\D^+_{\fg,W}(R)/\mathcal U_W$ be the Verdier quotient.
For $C^{\bullet}\in\D^+_{\fg,W}(R)\cap\D^{\geq a}(R)$, set
\[
 Z_n^{\bullet}=\tau_{\geq n}C^{\bullet}\quad(n\geq a),\qquad E^{\bullet}=\bigoplus_{n\geq a}Z_n^{\bullet},
\]
with coproduct injections $j_n$. The following are equivalent:
\begin{enumerate}[label=(\roman*),leftmargin=*]
\item $C^{\bullet}\in\mathcal U_W$;
\item $(qE^{\bullet},(qj_n)_n)$ is a coproduct of the family $(qZ_n^{\bullet})_n$;
\item the restriction map
\begin{equation}\label{eq:single-target-restriction}
 \Hom_{\mathcal Q}(qE^{\bullet},qC^{\bullet})\longrightarrow
 \prod_{n\geq a}\Hom_{\mathcal Q}(qZ_n^{\bullet},qC^{\bullet})
\end{equation}
is bijective.
\end{enumerate}
\end{proposition}
\begin{proof}
The family $(Z_n^{\bullet})$ is admissible. Let
$t_n:Z_n^{\bullet}\to Z_{n+1}^{\bullet}$ be the canonical quotient
maps, and define $\sigma:E^{\bullet}\to E^{\bullet}$ by
$\sigma j_n=j_{n+1}t_n$. In each degree, only finitely many tails
contribute to cohomology:
\[
 H^i(E^{\bullet})
 =\bigoplus_{n=a}^{i}H^i(C^{\bullet}),
\]
where the sum is zero when $i<a$. Under this decomposition,
$H^i(\sigma)$ shifts the summand indexed by $n$ to the summand
indexed by $n+1$ and kills the last summand. Thus
\[
 H^i(\sigma)^{\,i-a+1}=0\quad(i\geq a),\qquad
 (1-H^i(\sigma))^{-1}
 =\sum_{j=0}^{i-a}H^i(\sigma)^j.
\]
Thus $1-\sigma$ is a quasi-isomorphism.

The cone of the canonical map $u_n:C^{\bullet}\to Z_n^{\bullet}$
is bounded and supported in $W$. Hence $q(u_n)$ is an isomorphism.
Write $v_n:qZ_n^{\bullet}\to qC^{\bullet}$ for its inverse. Then
$v_{n+1}q(t_n)=v_n$.
Assume (iii). Surjectivity in \eqref{eq:single-target-restriction}
gives $l:qE^{\bullet}\to qC^{\bullet}$ with $lq(j_n)=v_n$.
The compatibility of the $v_n$ gives $l(1-q\sigma)q(j_n)=0$
for every $n$. Injectivity in \eqref{eq:single-target-restriction}
therefore gives $l(1-q\sigma)=0$.
Since $1-q\sigma$ is an isomorphism, $l=0$, and hence $qC^{\bullet}=0$.
The kernel of $q$ is $\mathcal U_W$, so (i) follows.
Conversely, (i) implies that every $Z_n^{\bullet}$ and their
admissible coproduct $E^{\bullet}$ belong to $\mathcal U_W$.
Their images are zero, giving (ii). The implication
(ii)$\Rightarrow$(iii) is the universal property of a coproduct.
\end{proof}

\section*{Acknowledgements}
The results and original mathematical proofs are the author's work.
ChatGPT Plus (model label: Sol) was used to find and check references,
check intermediate arguments, and revise the exposition and formulas.
The author is responsible for the mathematical claims and the final
manuscript.

\end{document}